\documentclass[a4paper]{amsart}
\RequirePackage{amsmath}
\RequirePackage{bm}
\RequirePackage{amssymb}
\RequirePackage{upref}
\RequirePackage{amsthm}
\RequirePackage{enumerate}
\RequirePackage{pb-diagram}
\RequirePackage{amsfonts}
\RequirePackage[mathscr]{eucal}
\RequirePackage{verbatim}
\RequirePackage{xr}
\RequirePackage{graphicx}
\usepackage{calc}
\usepackage{xspace}
\RequirePackage{color}
\RequirePackage{ifthen}
\RequirePackage{fancybox}
\RequirePackage{mathrsfs}

\newcommand{\cf}{cf.\@\xspace}
\newcommand{\resp}{resp.\@\xspace}

\newcommand{\al}{\alpha}
\newcommand{\bet}{\beta}
\newcommand{\ga}{\gamma}
\newcommand{\de}{\delta }
\newcommand{\e}{\epsilon}

\newcommand{\f}{\varphi}
\newcommand{\h}{\eta}

\newcommand{\ka}{\kappa}
\newcommand{\lam}{\lambda}
\newcommand{\m}{\mu}
\newcommand{\n}{\nu}
\newcommand{\om}{\omega}

\newcommand{\vt}{\vartheta}

\newcommand{\s}{\sigma}
\newcommand{\x}{\xi}

\newcommand{\C}{\varGamma}
\newcommand{\D}{\varDelta}

\newcommand{\Lam}{\varLambda}
\newcommand{\Om}{\varOmega}

\newcommand{\di}[1]{#1\nobreakdash-\hspace{0pt}dimensional}%\di n

\newcommand{\fv}[2]{#1\hspace{0pt}_{|_{#2}}}

\newcommand{\so}{{\mc S_0}}

\newcommand{\const}{\tup{const}}

\newcommand{\msp[1]}[1]{\mspace{#1mu}}

\newcommand{\R}[1][n+1]{{\protect\mathbb R}^{#1}}

\newcommand{\Hh}[1][n+1]{{\protect\mathbb H}^{#1}}

\newcommand{\N}{{\protect\mathbb N}}

\newcommand{\Z}{{\protect\mathbb Z}}
\newcommand{\eR}{\stackrel{\lower1ex \hbox{\rule{6.5pt}{0.5pt}}}{\msp[3]\R[]}}
\newcommand{\eN}{\stackrel{\lower1ex \hbox{\rule{6.5pt}{0.5pt}}}{\msp[1]\N}}
\newcommand{\eO}{\stackrel{\lower1ex \hbox{\rule{6pt}{0.5pt}}}{\msc O}}

\DeclareMathOperator{\graph}{graph}

\DeclareMathOperator{\B}{B}

\newcommand\im{\implies}
\newcommand\ra{\rightarrow}

\newcommand\hra{\hookrightarrow}

\newcommand\pa{\partial}
\newcommand\pde[2]{\frac {\partial#1}{\partial#2}}
\newcommand\pd[3]{\frac {\partial#1}{\partial#2^#3}}   %e.g. \pd fxi
\newcommand{\un}{\infty}
\newcommand{\A}{\forall}
\newcommand{\E}{\exists}

\newcommand{\set}[2]{\{\,#1\colon #2\,\}}
\newcommand{\uu}{\cup}
\newcommand{\ii}{\cap}
\newcommand{\uuu}{\bigcup}

\newcommand{\uud}{ \stackrel{\lower 1ex \hbox {.}}{\uu}}
\newcommand{\uuud}[1]{ \stackrel{\lower 1ex \hbox {.}}{\uuu_{#1}}}
\newcommand\su{\subset}

\newcommand\eS{\emptyset}
\newcommand{\sminus}[1][28]{\raise 0.#1ex\hbox{$\scriptstyle\setminus$}}

\newcommand\inn[1]{{\stackrel{\msp[9]\circ}{#1}}}

\newcommand{\wed}{\wedge}

\newcommand{\abs}[1]{\lvert#1\rvert}

\newcommand{\norm}[1]{\lVert#1\rVert}

\newcommand{\nnorm}[1]{| \mspace{-2mu} |\mspace{-2mu}|#1| \mspace{-2mu}
|\mspace{-2mu}|}
\newcommand{\spd}[2]{\protect\langle #1,#2\protect\rangle}

\newcommand\ch[3]{\varGamma_{#1#2}^#3}
\newcommand\cha[3]{{\bar\varGamma}_{#1#2}^#3}

\newcommand{\riem}[4]{R_{#1#2#3#4}}
\newcommand{\riema}[4]{{\bar R}_{#1#2#3#4}}

\newcommand{\tbf}{\textbf}
\newcommand{\tit}{\textit}

\newcommand{\tup}{\textup}% text upright

\newcommand{\mc}{\protect\mathcal}
\newcommand{\msc}{\protect\mathscr}

\newcommand\dt[1]{\dot{\tilde{#1}}}
\newcommand{\tlam}{\tilde\lam}

\providecommand{\bysame}{\makebox[3em]{\hrulefill}\thinspace}

\newcommand{\cq}[1]{\glqq{#1}\grqq\,}

\newcommand{\bt}{\begin{thm}}
\newcommand{\bl}{\begin{lem}}
\newcommand{\bc}{\begin{cor}}
\newcommand{\bd}{\begin{definition}}
\newcommand{\bpp}{\begin{prop}}
\newcommand{\br}{\begin{rem}}
\newcommand{\bn}{\begin{note}}
\newcommand{\be}{\begin{ex}}
\newcommand{\bes}{\begin{exs}}
\newcommand{\bb}{\begin{example}}
\newcommand{\bbs}{\begin{examples}}
\newcommand{\ba}{\begin{axiom}}
\newcommand{\bas}{\begin{assumption}}

\newcommand{\et}{\end{thm}}
\newcommand{\el}{\end{lem}}
\newcommand{\ec}{\end{cor}}
\newcommand{\ed}{\end{definition}}
\newcommand{\epp}{\end{prop}}
\newcommand{\er}{\end{rem}}
\newcommand{\en}{\end{note}}
\newcommand{\ee}{\end{ex}}
\newcommand{\ees}{\end{exs}}
\newcommand{\eb}{\end{example}}
\newcommand{\ebs}{\end{examples}}
\newcommand{\ea}{\end{axiom}}
\newcommand{\eas}{\end{assumption}}

\newcommand{\bp}{\begin{proof}}
\newcommand{\ep}{\end{proof}}
\newcommand{\eps}{\renewcommand{\qed}{}\end{proof}}

\newcommand{\bal}{\begin{align}}

\newcommand{\bi}[1][1.]{\begin{enumerate}[\upshape #1]}
\newcommand{\bia}[1][(1)]{\begin{enumerate}[\upshape #1]}
\newcommand{\bin}[1][1]{\begin{enumerate}[\upshape\bfseries #1]}
\newcommand{\bir}[1][(i)]{\begin{enumerate}[\upshape #1]}
\newcommand{\bic}[1][(i)]{\begin{enumerate}[\upshape\hspace{2\cma}#1]}
\newcommand{\bis}[2][1.]{\begin{enumerate}[\upshape\hspace{#2\parindent}#1]}
\newcommand{\ei}{\end{enumerate}}

\newcommand\ndots{\raise 0.47ex \hbox {,}\hskip0.06em\cdots %
     \raise 0.47ex \hbox {,}\hskip0.06em} 

\newcommand{\q}{\quad}
\newcommand{\qq}{\qquad}

\newcommand{\hp}{\hphantom}

\newcommand\nd{\noindent}

\newskip\Csmallskipamount                                                
\Csmallskipamount=\smallskipamount
\newskip\Cmedskipamount
\Cmedskipamount=\medskipamount
\newskip\Cbigskipamount
\Cbigskipamount=\bigskipamount

\newcommand\cvs{\vspace\Csmallskipamount}   
\newcommand\cvm{\vspace\Cmedskipamount}

\newskip\csa
\csa=\smallskipamount

\newskip\cma
\cma=\medskipamount

\newskip\cba
\cba=\bigskipamount

\newdimen\spt
\newcommand\citem{\cvs\advance\itemno by
1{(\romannumeral\the\itemno})\hskip3pt}
\newcommand{\bitem}{\cvm\nd\advance\itemno by
1{\bf\the\itemno}\hspace{\cma}}

\newcount\itemno
\newcommand{\las}[1]{\label{S:#1}}

\newcommand{\lae}[1]{\label{E:#1}}
\newcommand{\lat}[1]{\label{T:#1}}
\newcommand{\lal}[1]{\label{L:#1}}
\newcommand{\lad}[1]{\label{D:#1}}
\newcommand{\lac}[1]{\label{C:#1}}

\newcommand{\lap}[1]{\label{P:#1}}
\newcommand{\lar}[1]{\label{R:#1}}

\newcommand{\laas}[1]{\label{Ass:#1}}

\newcommand{\rt}[1]{Theorem~\ref{T:#1}}
\newcommand{\rl}[1]{Lemma~\ref{L:#1}}

\newcommand{\ras}[1]{Assumption~\ref{Ass:#1}}
\newcommand{\rc}[1]{Corollary~\ref{C:#1}}

\newcommand{\rp}[1]{Proposition~\ref{P:#1}}

\newcommand{\re}[1]{\eqref{E:#1}}
\newcommand{\frc}[1]{Corollary~\ref{C:#1} on page~\tup{\pageref{C:#1}}}
\newcommand{\frt}[1]{Theorem~\ref{T:#1} on page~\tup{\pageref{T:#1}}}
\newcommand{\frl}[1]{Lemma~\ref{L:#1} on page~\tup{\pageref{L:#1}}}
\newcommand{\frr}[1]{Remark~\ref{R:#1} on page~\tup{\pageref{R:#1}}}
\newcommand{\frd}[1]{Definition~\ref{D:#1} on page~\tup{\pageref{D:#1}}}
\newcommand{\fre}[1]{\eqref{E:#1} on page~\tup{\pageref{E:#1}}}

\newcommand{\frs}[1]{Section~\ref{S:#1} on page~\tup{\pageref{S:#1}}}

\newcommand{\fras}[1]{Assumption~\ref{Ass:#1} on page~\tup{\pageref{Ass:#1}}}

\newskip\thmskip
\thmskip=\parindent

\newskip\hsk
\newenvironment{hinw}{\labelsep=0pt\begin{list}{}{\labelsep=0pt\itemindent=0pt\labelwidth=0pt\leftmargin=\parindent\rightmargin=0pt\partopsep=\cba}%
\item\it\nopagebreak\nopagebreak}%
{\end{list}}

\newcommand\bh{\begin{hinw}}
\newcommand{\eh}{\end{hinw}}

\newtheoremstyle{normal}% name
  {\cba}%      Space above, empty = `usual value'
  {\cba}%      Space below
  {}% Body font
  {\thmskip}%Indent amount (empty = no indent, \parindent = para indent)
  {\bfseries}% Thm head font
  {.}%        Punctuation after thm head
  {\hsk}%     Space after thm head: " " = normal interword space;
  {}% Thm head spec

\newtheoremstyle{abschnitt}% name
  {\cba}%      Space above, empty = `usual value'
  {\cba}%      Space below
  {}% Body font
  {\thmskip}% Indent amount (empty = no indent, \parindent = para indent)
  {\bfseries}% Thm head font
  {.}%        Punctuation after thm head
  {\hsk}%     Space after thm head: " " = normal interword space;
  {}% Thm head spec

\newtheoremstyle{italic}% name
  {\cba}%      Space above, empty = `usual value'
  {\cba}%      Space below
  {\itshape}% Body font
  {\thmskip}%  Indent amount (empty = no indent, \parindent = para indent)
  {\bfseries}% Thm head font
  {.}%        Punctuation after thm head
  {\hsk}%     Space after thm head: " " = normal interword space;
  {}% Thm head spec

\newtheoremstyle{aufgaben}% name
  {\cba}%      Space above, empty = `usual value'
  {\cba}%      Space below
  {}% Body font
  {}%         Indent amount (empty = no indent, \parindent = para indent)
  {\normalsize\bfseries}% Thm head font
  {.}%        Punctuation after thm head
  {\hsk}%     Space after thm head: " " = normal interword space;
  {}% Thm head spec

\newtheoremstyle{break}% name
  {\cba}%      Space above, empty = `usual value'
  {\cba}%      Space below
  {\itshape}% Body font
  {}%         Indent amount (empty = no indent, \parindent = para indent)
  {\bfseries}% Thm head font
  {.}%        Punctuation after thm head
  {\newline}% Space after thm head: \newline = linebreak
  {}%         Thm head spec

\theoremstyle{italic}
\newtheorem{thm}[subsection]{Theorem}
\newtheorem{lem}[subsection]{Lemma}
\newtheorem{prop}[subsection]{Proposition}
\newtheorem{cor}[subsection]{Corollary}

\theoremstyle{normal}
\newtheorem{rem}[subsection]{Remark}
\newtheorem{definition}[subsection]{Definition}
\newtheorem{example}[subsection]{Example}
\newtheorem{examples}[subsection]{Examples}
\newtheorem{ex}[subsection]{Exercise}
\newtheorem{note}[subsection]{}
\newtheorem{axiom}[subsection]{Axiom}
\newtheorem{assumption}[subsection]{Assumption}

\theoremstyle{aufgaben}
\newtheorem{exs}[subsection]{Exercises}

\numberwithin{equation}{section}
\numberwithin{figure}{section}

\newenvironment{textequation}[1][0.8]
{\begin{equation}
\begin{aligned}
\begin{minipage}{#1\linewidth}}
{\end{minipage}
\end{aligned}
\end{equation}
\ignorespacesafterend}

\newcommand{\btext}{\begin{textequation}}
\newcommand{\etext}{\end{textequation}}

\def\hinweis{\@startsection{subsection}{2}%
 \z@{0.7\linespacing\@plus 0.5\linespacing}{0.7\linespacing}%
{\normalfont\itshape\indent}}

\newcounter{hours}\newcounter{minutes}
\newcommand{\printtime}{%
\setcounter{hours}{\time/60}%
\setcounter{minutes}{\time-\value{hours}*60}%
\ifthenelse{\value{minutes}<10}{\thehours :0\theminutes}{\thehours:\theminutes}}

\usepackage[german,english]{babel}
\usepackage{graphicx}
\RequirePackage{amsmath}
\RequirePackage{bm}
\RequirePackage{amssymb}
\RequirePackage{upref}
\RequirePackage{amsthm}
\RequirePackage{enumerate}
\RequirePackage{pb-diagram}
\RequirePackage{amsfonts}
\RequirePackage[mathscr]{eucal}
\RequirePackage{verbatim}
\RequirePackage{xr}
\RequirePackage{graphicx}
\RequirePackage{mathrsfs}
\usepackage{calc}
\usepackage{xspace}
\usepackage{dsfont}

\makeatletter
\RequirePackage{color}
\newcommand{\ann}[1]{\renewcommand{\@makefnmark}{\mbox{$^{\color{red}{\@thefnmark}}$}}%
\footnote {#1}}
\makeatother

\RequirePackage{upref}
\RequirePackage{amsthm}
\RequirePackage{enumerate}%\begin{enumerate}[(i)]
\usepackage[mathscr]{eucal}
\usepackage{xr-hyper}

\usepackage{calc}

\newlength{\oddsidemarginlength}
\newlength{\topmarginlength}

\newcounter{numberoflines}
\newcounter{tempcc}
\oddsidemargin=\oddsidemarginlength
\evensidemargin=\oddsidemargin
\topmargin=\topmarginlength
\usepackage[colorlinks=true,linkcolor=blue,citecolor=blue,urlcolor=blue]{hyperref}  

\begin{document}

\flushbottom

%\larger[1]
%\frontmatter

\title[CMC foliations of open spacetimes]{CMC foliations of  open spacetimes asymptotic to open Robertson-Walker spacetimes}

% author one information
\author{Claus Gerhardt}
\address{Ruprecht-Karls-Universit\"at, Institut f\"ur Angewandte Mathematik,
Im Neuenheimer Feld 205, 69120 Heidelberg, Germany}
%\curraddr{}
\email{\href{mailto:gerhardt@math.uni-heidelberg.de}{gerhardt@math.uni-heidelberg.de}}
\urladdr{\href{http://www.math.uni-heidelberg.de/studinfo/gerhardt/}{http://www.math.uni-heidelberg.de/studinfo/gerhardt/}}
%\thanks{This work was supported by the DFG}

% author two information
%\author{}
%\address{}
%\curraddr{}
%\email{}
%\thanks{}
%
%\subjclass[2000]{35J60, 53C21, 53C44, 53C50, 58J05}
%\keywords{Lorentzian manifold, mass, cosmological spacetime, general relativity, inverse mean curvature flow, ARW spacetimes}

\subjclass[2000]{35J60, 53C21, 53C44, 53C50, 58J05}
\keywords{Lorentzian manifold, open spacetimes, CMC foliation, general
relativity}\date{\today}
%
% at present the "communicated by" line appears only in ERA and PROC
%\commby{}

%\dedicatory{}

\begin{abstract} 
We consider open globally hyperbolic spacetimes $N$ of dimension $n+1$, $n\ge 3$, which are spatially asymptotic to a Robertson-Walker spacetime or  an open Friedmann universe with spatial curvature $\tilde\ka=0,-1$ and prove, under reasonable assumptions, that there exists a unique foliation by spacelike hypersurfaces of constant mean curvature and that the mean curvature function $\tau$ is a smooth time function if $N$ is smooth. Moreover, among the Friedmann universes which satisfy the necessary conditions are those that reflect the present assumptions of the development of the universe.
\end{abstract}

\maketitle

\tableofcontents

\setcounter{section}{0}
\section{Introduction}
Foliating a Lorentzian manifolds $N=N^{n+1}$ by spacelike hypersurfaces of constant mean curvature (CMC hypersurfaces) and using the mean curvature $\tau$ of the foliation hypersurfaces as a time function is very important for physical models of the universe. Solving the Einstein equations is a lot easier if the initial hypersurface has constant mean curvature and in Friedmann universes the mean curvature of the CMC hypersurfaces is also known as the Hubble constant---apart from a sign.

In case $N$ is globally hyperbolic and spatially compact, i.e., in case the Cauchy hypersurfaces are compact, the existence of a foliation by CMC hypersurfaces has been proved in \cite{cg1}. If the Cauchy hypersurfaces are non-compact only trivial CMC foliations in Robertson-Walker spacetimes are known so far. In this paper we prove the existence of  a CMC foliation in open globally hyperbolic spacetimes $N=N^{n+1}$, $n\ge 3$, which are spatially asymptotic to a Robertson-Walker spacetime $\hat N$. A Robertson-Walker spacetime is the warped product 
\begin{equation}
\hat N=I\times \hat\so,
\end{equation}
where $\hat\so$ is a space of constant curvature and we consider the cases where $\hat\so$ is either $\R[n]$ or $\Hh[n]$. These assumptions on the Cauchy hypersurfaces are also favoured in present cosmological models where it is mostly assumed that $\hat\so=\R[n]$. 

Let $\hat H=\hat H(t)$ be the mean curvature of the slices
\begin{equation}
\{x^0=t\},\qq\,t\in I,
\end{equation}
in $\hat N$ with respect to the past directed normal, i.e.,
\begin{equation}
\hat H=g^{ij}h_{ij}=-n \frac{a'}a,
\end{equation}
 where $a$ is the scale factor and $h_{ij}$ the second fundamental form, then the only condition we impose on $\hat N$ is 
\begin{equation}\lae{1.3.1} 
\hat H'>0,
\end{equation}
where a prime indicates differentiation with respect to $t$. This condition is satisfied by the models for an expanding Friedmann  universe, see e.g., \cite{cg:friedmann}, where the expansion is driven by dark matter and dark energy densities.

The existence of a CMC foliation is achieved by first solving the Dirichlet problems
\begin{equation}\lae{1.4}
\begin{aligned}
\fv{H}{M}&=\hat H(t_0)\qq \tup{in} \, B_R=B_R(\bar x_0),\\
\fv u{\pa B_R}&=t_0
\end{aligned}
\end{equation}
for spacelike graphs, $M=\graph u$, over nested balls  $B_R\su \so$ with fixed center and expanding radii, and proving uniform a priori estimates in $C^{m,\al}(\bar B_R)$, $m\ge 3$, $0<\al<1$, independent of $R$ and then letting $R$ tend to infinity. We then obtain a unique foliation of $N$ by spacelike hypersurfaces
\begin{equation}
M_t=\graph \fv{u((t,x))}{\so},\qq \,t\in I,
\end{equation}
having constant mean curvature $\hat H(t)$. The hypersurfaces $M_t$ uniformly converge to the slices
\begin{equation}
\{x^0=t\}
\end{equation}
in $C^{m,\al}(\so)$ if $r(x)$ tends to infinity, where $r(x)$ is a radial distance function. Finally, the mean curvature $\tau$ of the foliation hypersurfaces is a smooth global time function if $N$ is smooth. Here, is a more formal statement of this result:

\bt
The functions
\begin{equation}
u(t,x),\qq (t,x)\in I\times \so,
\end{equation}
describing the foliation by spacelike  hypersurfaces $M_t$, $t\in I$, are of class $C^{m-3,1}$ in $t$ such that
\begin{equation}
D^k_t u\in C^{m-k,\al}(\so)\qq\A\, 0\le k\le m-3,
\end{equation}
if $N$ is of class $C^{m,\al}$, $m\ge 3$, $0<\al<1$; if $N$ is smooth, i.e., $m=\un$, then $u$ is also smooth in the variables $(t,x)$ and the mean curvature function $\tau=\tau(x^0,x)$ is a smooth time function.
\et

\section{Notations, assumptions and definitions}\las{2} 

The main objective of this section is to state the equations of Gau{\ss}, Codazzi,
and Weingarten for spacelike hypersurfaces $M$ in a  \di{(n+1)} Lorentzian  space $N$.
Geometric quantities in $N$ will be 
denoted by
$(\bar g_{\al\bet}),(\riema \al\bet\ga\de)$, etc., and those in $M$ by $(g_{ij}), (\riem
ijkl)$, etc. Greek indices range from $0$ to $n$ and Latin from $1$ to $n$; the
summation convention is always used. Generic coordinate systems in $N$ resp.
$M$ will be denoted by $(x^\al)$ resp. $(\x^i)$. Covariant differentiation will
simply be indicated by indices, only in case of possible ambiguity they will be
preceded by a semicolon, i.e. for a function $u$ in $N$, $(u_\al)$ will be the
gradient and
$(u_{\al\bet})$ the Hessian, but e.g., the covariant derivative of the curvature
tensor will be abbreviated by $\riema \al\bet\ga{\de;\e}$. We also point out that
\begin{equation}
\riema \al\bet\ga{\de;i}=\riema \al\bet\ga{\de;\e}x_i^\e
\end{equation}
with obvious generalizations to other quantities.

Let $M$ be a \tit{spacelike} hypersurface, i.e. the induced metric is Riemannian,
with a differentiable unit normal $\n$ that is timelike.

In local coordinates, $(x^\al)$ and $(\x^i)$, the geometric quantities of the
spacelike hypersurface $M$ are connected through the following equations
\begin{equation}\lae{1.3}
x_{ij}^\al= h_{ij}\n^\al
\end{equation}
the so-called \tit{Gau{\ss} formula}. Here, and also in the sequel, a covariant
derivative is always a \tit{full} tensor, i.e.,
\begin{equation}
x_{ij}^\al=x_{,ij}^\al-\ch ijk x_k^\al+\cha \bet\ga\al x_i^\bet x_j^\ga.
\end{equation}
The comma indicates ordinary partial derivatives.

In this implicit definition the \tit{second fundamental form} $(h_{ij})$ is taken
with respect to $\n$.

The second equation is the \tit{Weingarten equation}
\begin{equation}
\n_i^\al=h_i^k x_k^\al,
\end{equation}
where we remember that $\n_i^\al$ is a full tensor.

Finally, we have the \tit{Codazzi equation}
\begin{equation}
h_{ij;k}-h_{ik;j}=\riema\al\bet\ga\de\n^\al x_i^\bet x_j^\ga x_k^\de
\end{equation}
and the \tit{Gau{\ss} equation}
\begin{equation}
\riem ijkl=- \{h_{ik}h_{jl}-h_{il}h_{jk}\} + \riema \al\bet\ga\de x_i^\al x_j^\bet x_k^\ga
x_l^\de.
\end{equation}

Now, let us assume that $N$ is a globally hyperbolic Lorentzian manifold with a
\tit{compact} Cauchy surface. $N$ is then a topological product $\R[]\times \mc
S_0$, where $\mc S_0$ is a compact Riemannian manifold, and there exists a
Gaussian coordinate system
$(x^\al)$, such that the metric in $N$ has the form
\begin{equation}\lae{1.7}
d\bar s_N^2=e^{2\psi}\{-{(dx^0)}^2+\s_{ij}(x^0,x)dx^idx^j\},
\end{equation}
where $\s_{ij}$ is a Riemannian metric, $\psi$ a function on $N$, and $x$ an
abbreviation for the spacelike components $(x^i)$, \cf \cite[Theorem 1.1]{sanchez:splitting}.
 We also assume that
the coordinate system is \tit{future oriented}, i.e. the time coordinate $x^0$
increases on future directed curves. Hence, the \tit{contravariant} timelike
vector$(\x^\al)=(1,0,\dotsc,0)$ is future directed as is its \tit{covariant} version
$(\x_\al)=e^{2\psi}(-1,0,\dotsc,0)$.

Let $M=\graph \fv u\so$ be the hypersurface
\begin{equation}
M=\set{(x^0,x)}{x^0=u(x),\,x\in \Om\su \mc S_0},
\end{equation}
where $\Om$ is an open domain in $\so$, then the induced metric has the form
\begin{equation}
g_{ij}=e^{2\psi}\{-u_iu_j+\s_{ij}\}
\end{equation}
where $\s_{ij}$ is evaluated at $(u,x)$, and its inverse $(g^{ij})=(g_{ij})^{-1}$ can
be expressed as
\begin{equation}\lae{1.10}
g^{ij}=e^{-2\psi}\{\s^{ij}+\frac{u^i}{v}\frac{u^j}{v}\},
\end{equation}
where $(\s^{ij})=(\s_{ij})^{-1}$ and
\begin{equation}\lae{1.11}
\begin{aligned}
u^i&=\s^{ij}u_j\\
v^2&=1-\s^{ij}u_iu_j\equiv 1-\abs{Du}^2.
\end{aligned}
\end{equation}
Hence, $\graph u$ is spacelike if and only if $\abs{Du}<1$.

The covariant form of a normal vector of a graph looks like
\begin{equation}
(\n_\al)=\pm v^{-1}e^{\psi}(1, -u_i).
\end{equation}
and the contravariant version is
\begin{equation}
(\n^\al)=\mp v^{-1}e^{-\psi}(1, u^i).
\end{equation}
Thus, we have
\br Let $M$ be spacelike graph in a future oriented coordinate system. Then, the
contravariant future directed normal vector has the form
\begin{equation}
(\n^\al)=v^{-1}e^{-\psi}(1, u^i)
\end{equation}
and the past directed
\begin{equation}\lae{1.15}
(\n^\al)=-v^{-1}e^{-\psi}(1, u^i).
\end{equation}
\er

In the Gau{\ss} formula \re{1.3} we are free to choose the future or past directed
normal, but we stipulate that we always use the past directed normal.

Look at the component $\al=0$ in \re{1.3} and obtain in view of \re{1.15}
\begin{equation}\lae{1.16}
e^{-\psi}v^{-1}h_{ij}=-u_{ij}-\cha 000\mspace{1mu}u_iu_j-\cha 0j0
\mspace{1mu}u_i-\cha 0i0\mspace{1mu}u_j-\cha ij0.
\end{equation}
Here, the covariant derivatives a taken with respect to the induced metric of
$M$, and
\begin{equation}
-\cha ij0=e^{-\psi}\bar h_{ij},
\end{equation}
where $(\bar h_{ij})$ is the second fundamental form of the hypersurfaces
$\{x^0=\const\}$.

An easy calculation shows
\begin{equation}
\bar h_{ij}e^{-\psi}=-\tfrac{1}{2}\dot\s_{ij} -\dot\psi\s_{ij},
\end{equation}
where the dot indicates differentiation with respect to $x^0$.

Sometimes, we need a Riemannian reference metric, e.g. if we want to estimate
tensors. Since the Lorentzian metric can be expressed as in \re{1.7}, we define a Riemannian reference metric $(\tilde g_{\al\bet})$ by
\begin{equation}
\tilde g_{\al\bet}dx^\al dx^\bet=e^{2\psi}\{{(dx^0)}^2+\s_{ij}dx^i dx^j\}
\end{equation}
and we abbreviate the corresponding norm of a vector field $\h$ by
\begin{equation}
\nnorm \h=(\tilde g_{\al\bet}\h^\al\h^\bet)^{1/2},
\end{equation}
with similar notations for higher order tensors.

Finally, let us recall the definition of an \tit{achronal} set:
\bd\lad{2.2}
A subset $M\su N$ is called \tit{achronal} provided any timelike curve meets $M$ at most once. For details see \cite[p. 413]{bn}.
\ed

Let us now formulate the assumptions on $N$. $N$ is a globally hyperbolic spacetime of dimension $n+1$ which is spatially asymptotic to a Robertson-Walker spacetime $\hat N$, which is  a warped product
\begin{equation}\lae{2.22}
\hat N=I\times \hat\so,
\end{equation}
where $\hat \so$ is either $\R[n]$ or $\Hh[n]$. The interval $I$ has the endpoints
\begin{equation}
I=(t_-,t_+).
\end{equation}
In physical applications
$\hat N$ is a Friedmann universe. In \cite{cg:friedmann} we proved the existence of an open Friedmann universe with
\begin{equation}
I=(0,\un)
\end{equation}
which has a big bang singularity and the mean curvature $\hat H(t)$ of the slices
\begin{equation}
\{x^0=t\}
\end{equation}
is negative with respect to the past directed normal such that
\begin{equation}\lae{2.26} 
\hat H'\ge c_0>0
\end{equation}
on compact subsets of $I$ and $\hat H$ satisfies
\begin{equation}
\lim_{t\ra 0}\hat H(t)\equiv H_-=-\un
\end{equation}
as well as
\begin{equation}
\lim_{t\ra\un}\hat H(t)\equiv H_+<0.
\end{equation}
Since
\begin{equation}
\hat H'=\abs{\hat A}^2+\hat R_{\al\bet}\nu^\al\nu^\bet=\abs{\hat A}^2+\hat R_{00},
\end{equation}
where
\begin{equation}
\abs{\hat A}^2=\hat h^{ij}\hat h_{ij},
\end{equation}
we have
\begin{equation}
\abs{\hat A}^2+\hat R_{\al\bet}\nu^\al\nu^\bet>0,
\end{equation}
while
\begin{equation}
\hat R_{00}=-n\frac{\Ddot a}a<0
\end{equation}
because 
\begin{equation}
\Ddot a>0,
\end{equation}
where $a(t)$ is the scale factor. 
\bas[Assumptions on $\hat N$]\laas{2.2}
In this paper we do not assume that $\hat N$ is a Friedmann universe nor do we assume that $\hat H$ is negative but we require \re{2.26} on compact subsets of $I$ or equivalently \re{1.3.1}. Then the limits 
\begin{equation}
\lim_{t\ra t_-}\hat H(t)\equiv H_-\ge-\un 
\end{equation}
and
\begin{equation}
\lim_{t\ra t_+}\hat H(t)\equiv H_+\le \un
\end{equation}
exist, i.e., $\hat N$ has only to satisfy \re{2.26}.
\eas
Since $N$ is supposed to be asymptotic to $\hat N$ we shall assume that both have the common time function $x^0$, that $N$ can also be written as a smooth product as in \re{2.22}, though we shall write $\so$ instead of $\hat\so$ because we consider $\so$ to be an embedded Cauchy hypersurface which is diffeomorphic to $\hat\so$. $N$ is supposed to be of class $C^{m,\al}$, $m\ge 3$, $0<\al<1$, i.e., the coordinate systems should be of class $C^{m,\al}$, the metric $\bar g_{\al\bet}$ of class $C^{m-1,\al}$ and the second fundamental form $\bar h_{ij}$ of the slices
\begin{equation}
\{x^0=\const\}
\end{equation}
of class $C^{m-2,\al}$. We assume that $N$ and $\hat N$ can be covered by a joint atlas of coordinate patches. The radial geodesic distance $\rho$ to a given point in $\hat\so$ is also defined in $\so$ but is of course only of class $C^{m,\al}$ there. In a coordinate slice
\begin{equation}
\{x^0=t_0\}
\end{equation}
in $\hat N$ the geodesic distance would be
\begin{equation}
r=a(t_0)\rho
\end{equation}
which is also defined when the corresponding slice is embedded in $N$.

The metric in $\hat N$ has the form
\begin{equation}
d\hat s^2=-(dx^0)^2+\hat\s_{ij}(x,x^0)dx^idx^j,
\end{equation}
where
\begin{equation}
\hat\s_{ij}=a^2(x^0)\tilde \s_{ij}(x)
\end{equation}
and $\tilde\s_{ij}$ is the metric in $\hat\so$.

The metric in $N$ can be written as
\begin{equation}
\begin{aligned}
d\bar s^2&=e^{2\psi}\{-(dx^0)^2+\s_{ij}(x,x^0)dx^idx^j\}\\
&\equiv\bar g_{\al\bet}dx^\al dx^\bet.
\end{aligned}
\end{equation}
Let  $\bar h_{ij}(x,x^0)$ \resp $\hat h_{ij}(x,x^0)$ be the second fundamental form of the slices
\begin{equation}
\{x^0=\const\}
\end{equation}
embedded in $N$ \resp $\hat N$, $e^{2\psi}\s_{ij}$ \resp $\hat \s_{ij}$ the corresponding induced metrics and $\bar \C^k_{ij}(x,x^0)$ \resp $\hat \C^k_{ij}(x,x^0)$ the corresponding Christoffel symbols, then we shall assume:
\bas[{Asymptotic behaviour}]\laas{2.3}
There exists $R_0>0$ and a constant $c>0$, which only depends on the compact sets in which $x^0$ ranges but not on $x$, provided 
\begin{equation}
r(x)>R_0,
\end{equation}
or equivalently,
\begin{equation}
\rho(x)>R_0,
\end{equation}
 such that
\btext
$\s_{ij}(x,x^0)$ and $\hat \s_{ij}(x,x^0)$ are uniformly equivalent as long as $x^0$ ranges in a compact subset of $I$.
\etext
\begin{equation}
\pm(\s_{ij}(x,x^0)-\hat\s_{ij}(x,x^0))\le \frac c{r(x)}\s_{ij}(x,x^0),
\end{equation}
\begin{equation}
\abs{e^\psi-1}\le \frac c{r(x)},
\end{equation}
\begin{equation}\lae{2.47}
\pm(\bar h_{ij}(x,x^0)-\hat h_{ij}(x,x^0))\le \frac c{r(x)}\s_{ij}(x,x^0),
\end{equation}
\begin{equation}\lae{2.48} 
\abs{\bar \C^k_{ij}(x,x^0)-\hat \C^k_{ij}(x,x^0)}\le \frac c{r(x)},
\end{equation}
where the norm on the left-hand side is an abbreviation for the norm of the corresponding tensor with respect to the metric $\s_{ij}(x,x^0)$. Furthermore, we assume
\begin{equation}\lae{2.49}
\abs{D^\ga(e^\psi-1)}\le \frac c{r(x)}\qq\A\, 0\le \abs\ga\le m-1,
\end{equation}
\begin{equation}
\abs{D^\ga(\bar h_{ij}(x,x^0)-\hat h_{ij}(x,x^0)}\le \frac c{r(x)}\qq\A\, 0\le \abs\ga\le m-2,
\end{equation}
and
\begin{equation}
\abs {D^\ga (\bar \C^k_{ij}-\hat \C^k_{ij})}\le \frac c{r(x)}\qq\A\, 0\le \abs\ga\le m-2,
\end{equation}
where the derivatives are either partial derivatives with respect to $x^0$ or spatial covariant derivatives with respect to $\s_{ij}(x,x^0)$.
\eas
\br
The previous assumptions on the asymptotic behaviour of $N$ and the assumption \re{2.26} imply
\begin{equation}\lae{2.52}
\dot{\bar H}(x,x^0)\ge c_0>0 \qq\A\,x\in\{r(x)>R_0\}
\end{equation}
uniformly in $x^0$ as long as $x^0$ stays in a compact subset of $I$.
\er

\bas[Additional assumptions on $N$]\laas{2.5}
$N$ should satisfy two additional assumptions. First, for any spacelike hypersurface of class $C^2$ we assume that
\begin{equation}\lae{2.53}
\abs A^2+\bar R_{\al\bet}\nu^\al\nu^\bet\ge0.
\end{equation}
Secondly, we assume the existence of future and past mean curvature barriers 
\begin{equation}
M_k=\graph \fv{\psi_k}\so,\qq k\in\Z,
\end{equation}
with mean curvatures
\begin{equation}
H_k=\fv H{M_k},
\end{equation}
which are, for fixed $k$, uniformly spacelike hypersurfaces satisfying 
\begin{equation}
\psi_k\in C^{3,\al}(\so),\qq\A\, k\in\Z,
\end{equation}
\begin{equation}
H_-<\inf_\so \fv H{ M_k}\le\sup_\so\fv H{M_k}<H_+,\qq\A\, k\in\Z,
\end{equation}
\begin{equation}
t_-=\lim_{k\ra-\un}\sup_\so \psi_k<\lim_{k\ra\un}\inf_\so\psi_k=t_+,
\end{equation}
and
\begin{equation}
H_-=\lim_{k\ra -\un}\sup_{M_k}H_k<\lim_{k\ra\un}\inf_{M_k} H_k=H_+
\end{equation}
where $t_\pm$ are the endpoints of $I$ and $H_\pm$ are the limits in \ras{2.2}. The previous assumptions on the $\psi_k$ also imply that we may assume without loss of generality
\begin{equation}
k<l\im \sup_\so\psi_k<\inf_\so\psi_l
\end{equation}
and
\begin{equation}
\sup_{M_k}H_k<\inf_{M_l}H_l,
\end{equation}
otherwise we consider a subsequence.
\eas

\section{$C^0$-estimates}

In this section we want to derive a priori bounds for solutions of the Dirichlet problem
\begin{equation}\lae{3.1}
\begin{aligned}
\fv{H}{M}&=f(u,x)\qq \tup{in} \, B_R(\bar x_0),\\
\fv u{\pa B_R}&=t_0
\end{aligned}
\end{equation}
where
\begin{equation}
M=\graph u=\{(x^0,x):x^0=u(x),\q x\in  B_R(\bar x_0)\}
\end{equation}
and 
\begin{equation}
B_R\equiv B_R(\bar x_0)\su \so.
\end{equation}
$\so$ is a Cauchy hypersurface. We shall generally assume that $\so$ is a fixed coordinate slice
\begin{equation}
\so=\{x^0=\bar t\}
\end{equation}
endowed with the induced metric 
\begin{equation}
e^{2\psi(\bar t,\cdot)}\s_{ij}(\bar t,\cdot).
\end{equation}

As a preparation let us first prove some lemmata.  
\bl\lal{3.1}
Let $M_i$, $i=1,2$, be spacelike hypersurfaces of class $C^1$ which are graphs over a bounded open domain $\Om\su \so$ such that
\begin{equation}\lae{3.6}
\inf_{\pa\Om}u_2>\sup_{\pa\Om}u_1
\end{equation}
and suppose, furthermore, that there exists a broken future directed timelike curve $\ga=(\ga^\al)$ of class $C^1$ from $\bar M_2$ to $\bar M_1$, then the endpoints of $\ga$ must both lie in the interior of the $M_i$.
\el
\bp
We argue by contradiction. Let the curve be parameterized over the interval $[0,1]$ such that
\begin{equation}
\ga(0)=p_2=(u(x_2),x_2)\in\bar M_2.
\end{equation}
and
\begin{equation}
\ga(1)=p_1=(u_1(x_1),x_1)\in \bar M_1.
\end{equation}

(i) First, assume that
\begin{equation}
x_2\in\pa\Om,
\end{equation}
then we have
\begin{equation}
u_1(x_1)=\ga^0(1)>\ga^0(0)=u_2(x_2)
\end{equation}
and hence
\begin{equation}
\ga(1)\in M_1,
\end{equation}
in view of the assumption \re{3.6}. Let $\tau_0$ be the largest $\tau$ such that
\begin{equation}
(\ga^i(\tau))\in\pa\Om,
\end{equation}
then
\begin{equation}
0\le\tau_0<1
\end{equation}
and
\begin{equation}
\ga^0(\tau_0)>\sup_{\pa\Om}u_1,
\end{equation}
because of \re{3.6}, from which we infer that there exists a future directed timelike curve $\tilde\ga$ parameterized over the interval $[0,1]$ such that
\begin{equation}
\tilde\ga(0)\in M_1
\end{equation}
and 
\begin{equation}
\tilde\ga(1)=\ga (1)\in M_1,
\end{equation}
such that $\tilde\ga$ is completely contained in the open cylinder
\begin{equation}
Q=I\times\Om
\end{equation}
which will lead to a contradiction as we shall show in the lemma below.

$\tilde\ga$ can be defined as follows: $\tilde\ga(0)$ is connected to $\ga(\tau)$ for  some $\tau>\tau_0$ by the timelike curve
\begin{equation}
\tilde\ga(s)=(s\ga^0(\tau)+(1-s)u_1(x_0),x_0),
\end{equation}
where
\begin{equation}
x_0=(\ga^i(\tau))\in\Om,\qq\tau>\tau_0,
\end{equation}
and $\tau$ has to satisfy
\begin{equation}
\ga^0(\tau)>u_1(\ga^i(\tau))
\end{equation}
which is valid if 
\begin{equation}
0<\tau-\tau_0<\de
\end{equation}
and $\de$ small enough.

After this first segment $\tilde\ga$ is identical with $\ga$. A reparameterization of  the first segment by setting
\begin{equation}
\tilde s=\tau s
\end{equation}
then leads the final definition of $\tilde\ga$.

(ii) Next, let us suppose
\begin{equation}
x_1\in \pa\Om,
\end{equation}
then
\begin{equation}
\ga(0)\in M_2.
\end{equation}
Let
\begin{equation}
\tau_0\in (0,1]
\end{equation}
be the first $\tau$ such satisfying
\begin{equation}
(\ga^i(\tau))\in \pa\Om,
\end{equation}
then
\begin{equation}
\ga^0(\tau_0)\le \ga^0(\tau_1)<\inf_{\pa\Om}u_2
\end{equation}
and we conclude that there exists a future directed timelike curve $\tilde\ga$ from
\begin{equation}
\ga(0)=p_2\in M_2
\end{equation}
to a point
\begin{equation}
\tilde\ga(1)\in M_2
\end{equation}
which lies completely inside the cylinder $Q$; again a contradiction.

$\tilde\ga$ is similarly defined as before, only, that now its last segment has to be defined by
\begin{equation}
\tilde\ga(s)=(su_2(\ga^i(\tau))+(1-s)\ga^0(\tau),\ga^i(\tau))
\end{equation}
for some $\tau<\tau_0$, where $\tau$ has to satisfy
\begin{equation}
u_2(\ga^i(\tau))>\ga,
\end{equation}
which is valid if
\begin{equation}
0<\tau_0-\tau<\de
\end{equation}
and $\de$ is small enough.
\ep

\bl\lal{3.2}
Let $M$ be a spacelike graph over an open bounded domain $\Om\su \so$,
\begin{equation}
M=\graph u,
\end{equation}
and assume that
\begin{equation}
u\in C^1(\bar\Om).
\end{equation}
Then, $M$ is achronal in the open cylinder
\begin{equation}
Q=I\times \Om.
\end{equation}
\el
\bp
For the definition of achronal confer \frd{2.2}. We argue by contradiction. Let $\ga=(\ga^\al)$ be a possibly broken timelike $C^1$-curve with image
\begin{equation}
\C=\{\ga(\tau):0\le\tau\le 1\} 
\end{equation}
and endpoints in $M$. Moreover, suppose that
\begin{equation}
\C\su Q.
\end{equation}
Without loss of generality let us assume that $\ga$ is future directed. Since $M$ is a $C^1$-graph it has a continuous timelike normal vector $\nu$ which we assume to be future directed.
The open set
\begin{equation}
U^+=\{(x^0,x):x^0>u(x),\q x\in\Om\}
\end{equation}
lies in the future of $M$ and
\begin{equation}
U^-=\{(x^0,x):x^0<u(x),\q x\in\Om\}
\end{equation}
in the past of $M$.

Let $\tau_0$ be the first $\tau>0$ such that
\begin{equation}
\ga(\tau)\in M,
\end{equation}
then the curve
\begin{equation}
\C_0=\{\ga(\tau):0\le\tau\le\tau_0\}
\end{equation}
intersects $M$ exactly twice, namely, at its endpoints. Since $\ga$ is future directed we deduce that there exists $\de>0$ such that
\begin{equation}
\ga(\tau)\in U^+\qq\A\, 0<\tau<\de
\end{equation}
and
\begin{equation}
\ga(\tau)\in U^-\qq\A\, \tau_0-\de<\tau<\de
\end{equation}
and we conclude that
\begin{equation}
\ga^0(\tau)-u(\ga^i(\tau))>0\qq\A\, 0<\tau<\de
\end{equation}
and
\begin{equation}
\ga^0(\tau)-u(\ga^i(\tau))<0\qq\A\, \tau_0-\de<\tau<\tau_0,
\end{equation}
hence, there exists
\begin{equation}
\tau_1\in (0,\tau_0)
\end{equation}
such that
\begin{equation}
\ga^0(\tau_1)=u(\ga^i(\tau_1))
\end{equation}
contradicting the fact that only the endpoints are part of $M$.
\ep
We are now ready to prove the crucial comparison theorem:
\bt\lat{3.3}
Let $\Om\su\so$ be a bounded open domain and let
\begin{equation}
M_i=\graph u_i,\qq i=1,2,
\end{equation}
be spacelike graphs over $\Om$ satisfying
\begin{equation}
u_i\in C^2(\bar \Om),
\end{equation}
\begin{equation}\lae{3.50} 
\inf_\Om \fv H{M_2}>\sup_\Om\fv H{M_1},
\end{equation}
and
\begin{equation}
\inf_{\pa\Om}u_2>\sup_{\pa\Om}u_1,
\end{equation}
then
\begin{equation}\lae{3.52}
u_1(x)<u_2(x)\qq\A\, x\in\Om.
\end{equation}
\et
\bp
First, we observe that the weaker conclusion
\begin{equation}\lae{3.53}
u_1\le u_2
\end{equation}
is as good as the stricter inequality \re{3.52} because of the weak Harnack inequality. Hence, suppose that \re{3.53} is not valid so that
\begin{equation}
E=\{x\in\bar\Om:u_2(x)<u_1(x)\}\not=\eS.
\end{equation}
Then, there exist points $p_i\in \bar M_i$ such that
\begin{equation}
\begin{aligned}
0<d_0&=d(\bar M_2,\bar M_1)=d(p_2,p_1)\\
&=\sup\{d(p,q):(p,q)\in \bar M_2\times \bar M_1\},
\end{aligned}
\end{equation}
where $d$ is the Lorentzian distance function which is continuous in globally hyperbolic spacetimes, \cf \cite[Lemma 4.5, p. 140]{beem-ehrlich}. Let $\ga=(\ga^\al)$ be a maximal future directed  geodesic from $\bar M_2$ to $\bar M_1$ realizing the distance with endpoints $p_i\in \bar M_i$, $i=1,2$, parameterized by arc length. Then, we first observe
\begin{equation}
p_i\in M_i,\qq \, i=1,2,
\end{equation}
in view of \rl{3.1}. 

We are now able to argue as in the proof of a corresponding result in \cite[Lemma 4.7.1]{cg:cp} to conclude that
\begin{equation}
\fv H{M_1}(p_1)\ge \fv H{M_2}(p_2),
\end{equation}
because of the assumption \fre{2.53}, contradicting \re{3.50}.
\ep
As a corollary we obtain: 
\bc
Let $\Om\su\so$ be a bounded open domain and $M=\graph u$ be a spacelike hypersurface, where $u\in C^2(\bar\Om)$ is a solution of the Dirichlet problem in $\Om$
\begin{equation}
\begin{aligned}
\fv HM&=f(u,x),\\
\fv u{\pa\Om}&=\f,
\end{aligned}
\end{equation}
satisfying
\begin{equation}
H_-<\inf_\Om f(u,x)\le \sup_\Om f(u,x)<H_+
\end{equation}
and
\begin{equation}
t_-<\inf_{\pa\Om}\f\le\sup_{\pa\Om}\f<t_+,
\end{equation}
then $u$ is a priori bounded if the assumption  \fre{2.53} is satisfied.
\ec
\bp
The proof follows immediately by employing appropriate barriers.
\ep

\section{Gradient estimates}
Let 
\begin{equation}
M=\graph u
\end{equation}
be a solution of the Dirichlet problem \fre{3.1} of class $C^{3}(\bar B_R)$. The right-hand side $f=f(x^0,x)$ should be of class $C^1$
\begin{equation}
f\in C^1(I\times\bar B_R),
\end{equation}
though the gradient estimates will not depend on $\pd fx0$ if
\begin{equation}
\pd fx0 \le 0,
\end{equation}
which is important when the penalization method is applied to approximate solutions of variational inequalities. We shall employ this method in the next section when  the existence of solutions to the Dirichlet problem is established. The fact that we consider a geodesic ball is of no importance. We could have chosen any precompact domain
\begin{equation}
\Om\su\so
\end{equation}
with $\pa\Om$ of class $C^2$. However, it is important that we consider constant boundary values for otherwise gradient estimates up to the boundary would be more difficult and the method we employ would fail. In \cite[Prop. 3.2]{br:mean} it is proved that non-constant boundary values can be reduced to the constant case  if certain conditions are satisfied but these conditions cannot  be verified easily.

We are going to prove the following theorem:
\bt\lat{4.1}
Let $M=\graph u$ be a solution of the Dirichlet problem \re{3.1} and let $K\su N$ be a compact set such that $M\su K$, then
\begin{equation}
\tilde v=v^{-1}=(1-\abs{Du}^2)^{-1}
\end{equation}
is uniformly bounded in $\bar B_R$,
\begin{equation}\lae{4.6}
\tilde v\le c,
\end{equation}
where $c$ depends on $K$, $t_0$ and the values $\abs f$, $\nnorm{Df}$ and ambient curvature terms in $K$ as well as the supremum norm of the mean curvature $\hat H_R$ of $\pa B_R$. If $K$ is bounded by coordinate slices
\begin{equation} 
t_1\le x^0\le t_2\qq\A\, (x^0,x)\in K,
\end{equation}
then $c$ depends on
\begin{equation}
c=c(t_1,t_2, \abs f_{\msc S}, \nnorm {Df}_{\msc S}, \abs{\hat H_R}_\msc S)
\end{equation}
and ambient curvature terms in
\begin{equation}
\msc S=\msc S(t_1,t_2),
\end{equation}
where $\msc S$ is the region defined by
\begin{equation}
\msc S=\{t_1\le x^0\le t_2\}
\end{equation}
provided the indicated norms in $\msc  S$ are  bounded. Since this is true in our case we can state that the a priori estimate is independent of $R$.
\et
The proof of the theorem requires some preparatory steps.

\bl\lal{4.2}
Let $M=\graph u$ satisfy the equation \fre{3.1}, then $\tilde v$ satisfies the elliptic differential equation
\begin{equation}\lae{3.20}
\begin{aligned}
-\D\tilde v=&-\norm A^2\tilde v
-f\h_{\al\bet}\n^\al\n^\bet\\
&-2h^{ij} x_i^\al x_j^\bet \h_{\al\bet}-g^{ij}\h_{\al\bet\ga}x_i^\bet
x_j^\ga\n^\al\\
&-\bar R_{\al\bet}\n^\al x_k^\bet\h_\ga x_l^\ga g^{kl}-f_\bet x_i^\bet  \h_\al x_k^\al g^{ik},
\end{aligned}
\end{equation}
where $\h$ is the covariant vector field $(\h_\al)=e^{\psi}(-1,0,\dotsc,0)$ and the covariant derivatives are to be understood with respect to the  induced metric $g_{ij}$ on $M$.
\el

\bp
We have $\tilde v=\spd \h\n$. Let $(\x^i)$ be local coordinates for $M$.
Differentiating $\tilde v$ covariantly we deduce
\begin{equation}\lae{3.21}
\tilde v_i=\h_{\al\bet}x_i^\bet\n^\al+\h_\al\n_i^\al,
\end{equation}
\begin{equation}\lae{3.22}
\begin{aligned}
\tilde v_{ij}= &\msp[5]\h_{\al\bet\ga}x_i^\bet x_j^\ga\n^\al+\h_{\al\bet}x_{ij}^\bet\n^\al\\
&+\h_{\al\bet}x_i^\bet\n_j^\al+\h_{\al\bet}x_j^\bet\n_i^\al+\h_\al\n_{ij}^\al
\end{aligned}
\end{equation}
Using then the Gau{\ss} formula, the Weingarten equation and the Ricci identities we obtain the desired result.
\ep

\bl\lal{3.3}
Let $K\su N$ be compact und $M\su K$, then there is a constant $c=c(K)$ such that for any positive function $0<\e=\e(x)$ on
$B_R$  we have
%\begin{equation}
\begin{align}
\nnorm \n&\le c\tilde v,\\\lae{3.14}
g^{ij}&\le c\tilde v^2\s^{ij},\\
\intertext{and}\lae{3.15}
\abs{h^{ij}\h_{\al\bet}x_i^\al x_j^\bet}&\le \frac{\e}{2}\norm A^2\tilde
v+\frac{c}{2\e}\tilde v^3
\end{align}
%\end{equation}
where $(\h_\al)$ is the vector field in \rl{4.2} and where we employed the Riemannian reference metric $\tilde g$. 
\el

\bp
The first two estimates can be immediately verified. To prove \re{3.15} we
choose local coordinates $(\x^i)$ such that
\begin{equation}
h_{ij}=\ka_i\de_{ij},\qq g_{ij}=\de_{ij}
\end{equation}
and deduce
\begin{equation}
\begin{aligned}
\abs{h^{ij}\h_{\al\bet}x_i^\al x_j^\bet}&\le \sum_i\abs{\ka_i}\abs{\h_{\al\bet}x_i^\al
x_i^\bet}\\
&\le \frac{\e}{2}\norm A^2\tilde v+\frac{1}{2\e}\tilde v^{-1}\sum_i\abs{\h_{\al\bet}
x_i^\al x_i^\bet}^2,
\end{aligned}
\end{equation}
and
\begin{equation}
\sum_i\abs{\h_{\al\bet}
x_i^\al x_i^\bet}^2\le g^{ik}\h_{\al\bet}x_i^\al x_j^\bet \msp[3]g^{jl} \h_{\ga\de} x_k^\ga
x_l^\de.
\end{equation}
Hence, the result in view of \re{3.14}.
\ep

Combining the preceding lemmata we infer 
\bl\lal{3.4}
There is a constant $c=c(K)$ such that for any positive function $\e=\e(x)$ on
$B_R$ the term $\tilde v$ satisfies a parabolic inequality of the form
\begin{equation}\lae{4.20}
-\D\tilde v\le -(1-\e)\norm A^2\tilde v+c[\abs f+\nnorm{Df}]\tilde
v^2+c[1+\e^{-1}]\tilde v^3.
\end{equation}
\el

We note that the statement \tit{$c$ depends on $K$} also implies that $c$
depends on geometric quantities of the ambient space restricted to $K$.

We further need the following two lemmata

\bl\lal{3.5}
Let $M=\graph u$ have prescribed mean curvature $f$, then 
\begin{equation}
-\D u=e^{-\psi}v^{-1}f-e^{-\psi}g^{ij}\bar h_{ij}+\cha 000\norm{Du}^2+2\cha
0i0 u^i.
\end{equation}
\el

\bp
This follows immediately from equation \re{1.16}. 
\ep

\bl\lal{3.6}
Let $M\su K$ be a graph over $B_R$, $M=\graph u$, then
\begin{equation}
\abs{\tilde v_i u^i}\le c\tilde v^3+\norm A e^\psi\norm {Du}^2,
\end{equation}
where $c=c(K)$.
\el

\bp
First, we use that
\begin{equation}
\tilde v^2=1+e^{2\psi}\norm{Du}^2,
\end{equation}
and thus,
\begin{equation}
2\tilde v\tilde v_i=2\psi_\al x_i^\al e^{2\psi} \norm{Du}^2+2e^{2\psi} u_{ij}u^j,
\end{equation}
from which we infer
\begin{equation}
\abs{\tilde v_iu^i}\le c\tilde v^3+\tilde v^{-1}e^{2\psi}\abs{u_{ij}u^i u^j},
\end{equation}
which gives the result because of \re{1.16}.
\ep

We are now ready to prove \rt{4.1}.

\bp[\tbf{Proof of \rt{4.1}}]
Let $\m,\lam$ be positive constants, where $\m$ is supposed to be small and $\lam$
large, and define
\begin{equation}\lae{3.27}
\f=e^{\m e^{\lam (u+c)}},
\end{equation}
where $c$ is large positive constant such that $u+c>1$.

We shall show that
\begin{equation}
w=\tilde v \f
\end{equation}
is uniformly bounded if $\m,\lam$ are chosen appropriately. Let $p_0=(u(y_0),y_0)\in \bar M$ be point where $w$ attains its supremum
\begin{equation}
w(p_0)=\sup_{\bar M}w,
\end{equation}
where we now consider the functions to be defined on $M$. We shall apply the maximum principle, or, at the boundary, a slight modification of it,  to obtain an a priori bound for $w$. 

(i) Let us first assume that $p_0\in M$, then we can apply the maximum principle directly.
In view of \rl{3.3} and \rl{3.5} we have
\begin{equation}
-\D\f\le c\m\lam e^{\lam(u+c)}[\tilde v\abs f +\tilde v^2] \f-\m\lam^2 e^{\lam(u+c)} [1+\m
e^{\lam(u+c)}]\norm{Du}^2\f,
\end{equation}
from which we further deduce taking \rl{3.4} and \rl{3.6} into account 
\begin{equation}
\begin{aligned}
-\D w&\le -(1-\e) \norm A^2\tilde v\f +c[\abs f+\nnorm{Df}]\tilde v^2\f\\
&\q\,+c[1+\e^{-1}]\tilde v^3\f-\m\lam^2 e^{\lam(u+c)} [1+\m e^{\lam(u+c)}] \tilde v
\norm{Du}^2\f\\
&\q\,+c[1+\abs f]\m\lam e^{\lam(u+c)}\tilde v^3\f+2\m\lam e^{\lam(u+c)} \norm A e^\psi
\norm{Du}^2\f.
\end{aligned}
\end{equation}

We estimate the last term on the right-hand side by
\begin{equation}
\begin{aligned}
2\m\lam e^{\lam(u+c)}\norm A e^\psi\norm{Du}^2\f&\le (1-\e)\norm A^2\tilde v\f\\
&\q\,+\frac{1}{1-\e}\m^2\lam^2e^{2\lam(u+c)}\tilde v^{-1}e^{2\psi}\norm{Du}^4\f,
\end{aligned}
\end{equation}
and conclude
\begin{equation}
\begin{aligned}
-\D w&\le  c[\abs f+\nnorm{Df}]\tilde v^2\f+ c[1+\abs f]\m\lam e^{\lam(u+c)} \tilde
v^3\f\\
 &\q\,+c[1+\e^{-1}]\tilde v^3\f +[\frac{1}{1-\e}-1]\m^2\lam^2 e^{2\lam
(u+c)}\norm{Du}^2\tilde v\f\\
&\q\,-\m\lam^2 e^{\lam(u+c)}\norm{Du}^2\tilde v\f,
\end{aligned}
\end{equation}
where we have used that
\begin{equation}
e^{2\psi}\norm{Du}^2\le \tilde v^2.
\end{equation}

Setting $\e=e^{-\lam (u+c)}$, we then obtain
\begin{equation}\lae{3.34}
\begin{aligned}
-\D w&\le c[\abs f+\nnorm{Df}]\tilde v^2\f+c e^{\lam(u+c)} \tilde v^3\f\\
&\q\,+c[1+\abs f]\m\lam e^{\lam(u+c)}\tilde v^3\f\\
&\q\,+[\frac{\m}{1-\e}-1]\m\lam^2 e^{\lam(u+c)}\norm{Du}^2\tilde v\f.
\end{aligned}
\end{equation}

Now, we choose $\m=\frac{1}{2}$ and $\lam_0$ so large that
\begin{equation}
\frac{\m}{1-e^{-\lam(u+c)}}\le \frac{3}{4}\qq\A\,\lam\ge \lam_0,
\end{equation}
and infer that the last term on the right-hand side of \re{3.34} is less than
\begin{equation}
-\frac{1}{8}\lam^2e^{\lam(u+c)}\norm{Du}^2\tilde v\f
\end{equation}
which in turn can be estimated from above by
\begin{equation}
-c\lam^2e^{\lam(u+c)}\tilde v^3\f
\end{equation}
at points where $\tilde v\ge 2$.

Thus, we conclude that for
\begin{equation}
\lam\ge \max (\lam_0, 4[1+\abs f_{_K}])
\end{equation}
the maximum principle, applied to $w$ at $y_0$, yields
\begin{equation}\lae{4.39}
w\le \const (\lam_0,\abs f, \nnorm{Df},K).
\end{equation}

\br
The estimate  \re{4.39} is also valid if $\f$ is replaced by 
\begin{equation}
\tilde \f=e^{\mu e^{-\lam(u-c+2t_0)}}
\end{equation}
and then considering
\begin{equation}
\tilde w=\tilde v\tilde \f
\end{equation}
instead of $w$, where $\lam,\mu$ and $c$ are defined as before. If $\tilde w$ attains its supremum at an interior point $p_0\in M$, then the inequality \re{4.39} is also valid if $w$ is replaced by $\tilde w$ on the left-hand side.
\er

(ii) We now assume that the supremum of $w$ is attained on the boundary of $M$, where
\begin{equation}
u=t_0,
\end{equation}
i.e., $u$ is constant on $\pa M$. We then argue similarly as Bartnik in the proof of \cite[Theorem 3.1]{br:mean}. Consider $w$ and $\tilde w$ simultaneously. If one of the functions attains its maximum in the interior, then the estimate is already proved. Thus, let us assume that both attain their maximum on $\pa M$. Since
\begin{equation}
w=\tilde w\qq \tup{on} \, \pa  M
\end{equation}
we consider a point $p_0\in\pa M$, or equivalently, $y_0\in\pa B_R$, such that $w$ and $\tilde w$ both attain their maximum in $p_0$. We may also assume that 
\begin{equation}
Du(y_0)\not=0,
\end{equation}
for otherwise we have nothing to prove. Hence, either
\begin{equation}
\frac{Du}{\norm {Du}}
\end{equation}
or
\begin{equation}
-\frac{Du}{\norm {Du}}
\end{equation}
is equal to the outer normal $\tilde \nu$ of $\pa M$, i.e.,
\begin{equation}
\tilde\nu\in T_{p_0}(M),
\end{equation}
since $u$ is constant on $\pa M$.
Let us assume that
\begin{equation}
\tilde\nu=-\frac{Du}{\norm {Du}},
\end{equation}
then we have in $p_0$
\begin{equation}
\begin{aligned}
0\le w_i\tilde \nu^i=\tilde v_i\tilde\nu^i\f +\tilde v\lam e^{\lam(u+c)}\f u_i\tilde\nu^i
\end{aligned}
\end{equation}
from which we deduce
\begin{equation}
\begin{aligned}
\lam e^{\lam(u+c)}\tilde v\norm{Du}&\le \f\abs{\tilde v_iu^i}\norm{Du}^{-1}\\
&=\f \abs{\tilde v_i\check u^i}\tilde v^2e^{-2\psi}\norm{Du}^{-1},
\end{aligned}
\end{equation}
where $\check u^i$ indicates that the index is raised with respect to the metric $\s_{ij}(u,x)$
\begin{equation}
\check u^i=\s^{ij}u_j
\end{equation}
and where we recall that
\begin{equation}
u^i=\tilde v^2 e^{-2\psi}\check u^i,
\end{equation}
\cf \cite[equ. (5.4.20) \& (5.4.21)]{cg:cp}.

On the other hand,,
\begin{equation}\lae{4.53}
\abs{\tilde v_i\check u^i}\le c\abs f+c\abs{\bar A}\tilde v+c\nnorm{D\psi}+ c\abs{\msc H}\tilde v,
\end{equation} 
where the constant $c$ only depends on geometric quantities of the ambient metric in the capped region $\msc S(t_1,t_2)$ and where $\msc H$ is the mean curvature of the geodesic sphere $\pa B_R$ embedded in $(\so, \s_{ij}(t_0,x))$, $\bar A$ is the second fundamental form of the coordinate slice $M_{t_0}$. A proof is given in the lemma below.

Since $\abs{\msc H}, \abs{\bar A}, \abs f$ and $\abs\psi,\nnorm{D\psi}$ are uniformly bounded in $\msc S$ we conclude 
\begin{equation}\lae{4.54}
\abs{\tilde v_i\check u^i}\le c_0\tilde v
\end{equation}
if $\tilde v>2$ and hence we deduce
\begin{equation}
\tilde v(p_0)<2
\end{equation}
if $\lam$ is larger than some constant $\lam_0$, where
\begin{equation}
\lam_0=\lam_0(\abs{\msc H}, \abs{\bar A}, \abs f, \abs\psi,\nnorm{D\psi}, \abs u)
\end{equation}
and where the norms on the right-hand side are supremum norms over $\pa B_R$ \resp $\msc S$ or $B_R$. Hence, we finally proved the a priori estimate for \re{4.6} in view of the result in the lemma below. 
\ep

\bl\lal{4.8}
Assume that $w$ attains its maximum at a point $p_0=(u(y_0),y_0)$ with $y_0\in \pa B_R$ and that
\begin{equation}\lae{4.57}
-\frac{\check u^i}{\abs{Du}}=\hat\nu^i, 
\end{equation}
where $\hat\nu^i$ is the outward normal to $\pa B_R$. Suppose, furthermore, that $\tilde v\ge 2$, then the estimate \re{4.54} is valid.
\el
\bp
First, let $\tilde g_{\al\bet}$ be the conformal metric to $\bar g_{\al\bet}$ such that
\begin{equation}
\bar g_{\al\bet}=e^{2\psi}\tilde g_{\al\bet}
\end{equation}
and consider $M$ to be embedded in $N$ equipped with the metric $\tilde g_{\al\bet}$ instead of $\bar g_{\al\bet}$. Denote the corresponding geometric quantities of $M$ by $\tilde g_{ij}, \tilde h_{ij}, \tilde H$ and $\tilde \nu =(\tilde \nu^\al)$. Then, we have
\begin{equation}\lae{4.59}
e^\psi H=\tilde H +n\psi_\al \tilde\nu^\al,
\end{equation}
\cf \cite[equ. (1.1.52) on p. 7]{cg:cp}, and
\begin{equation}
g^{ij}=e^{-2\psi}\tilde g{ij}=e^{-2\psi}(\s^{ij}+\tilde v^2\check u^i\check u^j).
\end{equation}
Let $u_{ij}$ be the second covariant derivatives of $u$ with respect to the metric $\tilde g_{ij}$ and $u_{;ij}$ the covariant derivatives with respect to the metric $\s_{ij}(x,u)$, then
\begin{equation}\lae{4.61}
u_{ij}=\tilde v^2u_{;ij},
\end{equation}
\cf \cite[Lemma 2.7.6]{cg:cp}. 

Next, let us differentiate
\begin{equation}
\tilde v=(1-\abs{Du}^2)^{-1}
\end{equation}
covariantly with respect to the metric $\s_{ij}(x,u)$ yielding
\begin{equation}
\tilde v_i=\tilde v^3 u_{;ij}\check u^j
\end{equation}
and hence
\begin{equation}\lae{4.64}
\tilde v_i\check u^i=\tilde v^3 u_{;ij}\check u^i\check u^j=\tilde v^3 u_{;ij} \hat\nu^i\hat\nu^j\abs{Du}^2.
\end{equation}

From \re{4.61} we infer that 
\begin{equation}\lae{4.65}
\begin{aligned}
\tilde vH&=-\tilde v^2\tilde g^{ij}u_{;ij}+\tilde g^{ij}\bar h_{ij}\\
&=-\tilde v^2\{\s^{ij}u_{;ij}+u_{;ij}\check u^i\check u^j\tilde v^2\}+\tilde g^{ij}\bar h_{ij},
\end{aligned}
\end{equation}
where $\bar h_{ij}$ is the second fundamental form of the coordinate slices with respect to the ambient metric $\tilde g_{\al\bet}$. Choosing a coordinate system $(x^i)$ in $y_0\in \pa B_R$ such that
\begin{equation}\lae{4.66}
\begin{aligned}
\s_{ij}&=\de_{ij},\\
\pd{}xn&=\hat\nu,
\end{aligned}
\end{equation}
and $T_{y_0}(\pa B_R)$ is spanned by
\begin{equation}
\pd{}xi,\qq 1\le i\le n-1,
\end{equation}
we see that we only have to estimate the tangential second derivatives of $u$ appropriately to complete the proof of the lemma. 

Let $x^i(\xi^a)$, $1\le a\le n-1$, be a local embedding of $\pa B_R$ into $(\so,\s_{ij}(t_0,x))$ around $y_0=x(\xi_0)$, then
\begin{equation}
0=u_a=u_ix^i_a
\end{equation}
and
\begin{equation}\lae{4.69}
0=u_{ab}=u_{:ij}x^i_ax^j_b-\hat h_{ab}\hat u_i\nu^i,
\end{equation}
where the semicolon indicates covariant differentiation with respect to the metric $\s_{ij}(t_0,x)$ and where $\hat h_{ab}$ is the second fundamental form of $\pa B_R$ with respect to the outward normal. Note that
\begin{equation}
u_{;ij}=u_{:ij}+c_{ij}\abs{Du}^2,
\end{equation}
where $c_{ij}$ is a bounded tensor.  Choosing now the coordinates $(\xi^a)$ at $\xi_0$ such that
\begin{equation}\lae{4.71}
x^i_a=\de^i_a
\end{equation}
we conclude
\begin{equation}\lae{4.72}
-\sum_{i=1}^{n-1}\s^{ij}u_{;ij}=\msc H\abs{Du}^2+c\abs{Du}^2,
\end{equation}
in view of \re{4.57}, \re{4.66} and \re{4.71}.

Combining now the relations  \re{4.72} and \re{4.65} we obtain in $y_0$
\begin{equation}
u_{;ij}\hat \nu^i\hat\nu^j(1+\tilde v^2\abs{Du}^2)=-v\tilde H+v^2(\bar H+\bar h_{ij}\check u^i\check u^j\tilde v^2)+\msc H\abs{Du}^2+c\abs{Du}^2
\end{equation}
from which we immediately infer the estimate \re{4.53}, in view of  \re{4.59}, \re{4.64} and the assumption that 
\begin{equation}
\abs{Du}^2\ge \frac12.
\end{equation}
Let us also recall that 
\begin{equation}
v={\tilde v}^{-1}.
\end{equation}
\ep

\section{Existence of a solution of the Dirichlet problem}

in this section we want to prove that the Dirichlet problem \fre{3.1} has a solution of class $C^{3,\al}$ where we assume that the function $f=f(x^0,x)$ is of class $C_{\tup{loc}}^{1,\al}(N)$ and satisfies the estimate
\begin{equation}
H_-<\inf_N f\le\sup_N f<H_+.
\end{equation}
Let 
\begin{equation}\lae{5.2}
M_i=\graph\psi_i,\qq\, i=1,2,
\end{equation}
be barriers satisfying
\begin{equation}
\psi_1<\psi_2,
\end{equation}
\begin{equation}
H_2=\fv H{M_2}>\sup f,
\end{equation}
\begin{equation}
\inf \psi_2>t_0,
\end{equation}
\resp
\begin{equation}
H_1=\fv H{M_1}<\inf f,
\end{equation}
and
\begin{equation}
\sup\psi_1<t_0.
\end{equation}
Then, we pick two constants $m_i\in I$ such that
\begin{equation}
m_1<\inf\psi_1<\sup\psi_2<m_2
\end{equation}
and $\e>0$ so small such that
\begin{equation}
m_2+\e\in I
\end{equation}
and
\begin{equation}
m_1-\e\in I.
\end{equation}
Furthermore, let $\bet_i$ be smooth real functions satisfying
\begin{equation}
\bet_1(t)=
\begin{cases}
-1&, t\le -\e,\\
\hp-0&, t\ge 0,
\end{cases}
\end{equation}
and
\begin{equation}
\bet_2(t)=
\begin{cases}
0&, t\le 0,\\
1&, t\ge \e,
\end{cases}
\end{equation}
as well as
\begin{equation}
\bet_i'\ge 0.
\end{equation}
Then, we look at the Dirichlet problem for $M=\graph u$ in $B_R$
\begin{equation}\lae{5.14}
\begin{aligned}
\fv HM&=f(x,u)-\tilde\mu \bet_1(u-m_1)-\tilde\mu\bet_2(u-m_2),\\
\fv u{\pa B_R}&=t_0,
\end{aligned}
\end{equation}
where $\tilde\mu>0$ is a sufficiently large constant. The exact value will be determined later. It will depend on $f$ and on geometric quantities of the ambient space in the capped region
\begin{equation}
\msc S(m_1-\e,m_2+\e)=[m_1-\e,m_2+\e]\times\so.
\end{equation}
We shall prove that the Dirichlet problem \re{5.14} has a solution 
\begin{equation}
u\in C^{1,\al}(\bar B_R)
\end{equation}
satisfying the estimates
\begin{equation}
m_1-\e\le u\le m_2+\e
\end{equation}
and
\begin{equation}
\tilde v\le c,
\end{equation}
where $c$ depends on $\tilde\mu$ but not on $\e$ and $\bet_i'$.

Then, if we would $\e$ letting tend to zero, we would obtain a solution of a variational inequality, where the obstacles are given by the slices
\begin{equation}
\{x^0=m_i\}.
\end{equation}
However, we do not need to do this, since, by employing the barriers and the comparison theorem, \frt{3.3}, we shall show that the penalization terms $\bet_i$ vanish. But we shall later refer to the above argument, when we claim that a given variational inequality has a solution.

Assuming for the moment that the Dirichlet problem \re{5.14} has a solution we are then able to prove:
\bt\lat{5.1}
The solution of the Dirichlet problem \re{5.14} actually satisfies the equation
\begin{equation}\lae{5.20}
\fv HM=f(u,x),
\end{equation}
since $u$ lies between the barriers $\psi_i$ in \re{5.2}
\begin{equation}\lae{5.21}
m_1<\psi_1<u<\psi_2<m_2
\end{equation}
and hence the penalization terms vanish.
\et
\bp
We shall only prove the upper estimate, since the argument for the lower estimate is similar, or even identical by simply  switching the light-cone. Also note that proving the weaker inequality
\begin{equation}
\psi\le u\le \psi_2
\end{equation}
is as good as the stronger inequality, since the $\bet_i$ terms would still vanish and the maximum principle would yield the final result.

Thus, let us assume that
\begin{equation}
\{x\in B_R:\psi_2(x)<u(x)\}\not=\eS.
\end{equation}
From the proof of \frt{3.3} we then conclude that there exists a maximizing future directed geodesic $\ga=(\ga^\al)$ from a point $p_2\in M_2$ to a point $p_1\in M$. Let $p_1$ have the coordinates $(u(x_1),x_1)$ and  similarly express $p_2$  as $(\psi_2(x_2),x_2)$. Then we deduce
\begin{equation}
u(x_1)>\psi_2(x_2)>\psi_1(x_2)>m_1
\end{equation}
and this estimate is also valid in a small ball $B_\rho(x_1)$, hence we have
\begin{equation}
u-m_1>0\qq\tup{in}\, B_\rho(x_1)
\end{equation}
and, in view of the definition of $\bet_1$, we infer
\begin{equation}
\bet_1(u-m_1)=0,
\end{equation}
hence,
\begin{equation}
\fv HM = f(u,x)-\tilde\mu\bet_2(u-m_2)\le f(u,x)
\end{equation}
in $B_\rho(x_1)$ and the further arguments in the proof of \rt{3.3} lead to a contradiction.
\ep

Let us now prove that the Dirichlet problem \re{5.14} has a solution $u\in C^{3,\al}(\bar B_R)$. We shall argue similarly as in the proof of \cite[Theorem 5.2]{cg1}, where we treated a related problem. Since we then considered compact hypersurfaces without boundary and also used a different method for the gradient estimate, we cannot simply refer to that result but have to argue a little more detailed.

For the existence proof we shall use a Leray-Schauder type fixed point argument. For technical reasons it is therefore necessary to consider $M$ embedded in the conformal space $(\tilde N,\tilde g_{\al\bet})$, where $\bar g_{\al\bet}$ and $\tilde g_{\al\bet}$ are related by
\begin{equation}
\bar g_{\al\bet}=e^{2\psi}\tilde g_{\al\bet}.
\end{equation}
We already used this embedding of $M$ in the proof of \frl{4.8}. Now, we need a more detailed description how the geometric quantities of the two embeddings are related. As before we embellish the geometric symbols in $\tilde N$ by a tilde. The geometric quantities are related by
\begin{equation}
g_{ij}=e^{2\psi}\tilde g_{ij},
\end{equation}
\begin{equation}
\nu=e^{-\psi}\tilde\nu,
\end{equation}
\begin{equation}
h_{ij}e^{-\psi}=\tilde h_{ij}+\psi_\al \tilde\nu^\al\tilde g_{ij},
\end{equation}
\begin{equation}
e^\psi H=\tilde H+n\psi_\al\tilde\nu^\al,
\end{equation}
\cf \cite[Proposition 1.1.11]{cg:cp}.

Let $u_{ij}$ be the covariant second derivatives of $u$ with respect to $\tilde g_{ij}$, then
\begin{equation}
u_{ij}=\tilde v^2 u_{;ij},
\end{equation}
where the covariant derivatives on the right-hand side are with respect to the metric $\s_{ij}(u,x)$, \cf \cite[Lemma 2.7.6]{cg:cp}, and hence
\begin{equation}
\tilde h_{ij}\tilde v=-\tilde v^2u_{;ij}+\bar h_{ij},
\end{equation}
where $\bar h_{ij}$ is the second fundamental form of the coordinate slices
\begin{equation}
\{x^0=\const\}
\end{equation}
in $\tilde N$---here, we refrain to embellish the barred quantities by an additional tilde---and we further deduce
\begin{equation}\lae{5.36} 
\tilde H\tilde v=-\tilde v^2\tilde g^{ij}u_{;ij}+\bar H +\tilde v^2 \bar h_{ij} \check u^i\check u^j,
\end{equation}
where we recall
\begin{equation}
\check u^i=\s^{ij}u_j.
\end{equation}

Thus, the Dirichlet problem \re{5.14} is equivalent to 
\begin{equation}\lae{5.38}
\begin{aligned}
\fv{\tilde H}M&=e^{\psi}\{f(u,x)-\tilde\mu \bet_1(u-m_1)-\tilde\mu\bet_2(u-m_2)+nD_\al(e^{-\psi})\tilde\nu^\al\},\\
\fv u{\pa B_R}&=t_0.
\end{aligned}
\end{equation}
Let us also note that
\begin{equation}\lae{5.39} 
\tilde v \tilde g^{ij}u_{;ij}=D_i(\tilde v\s^{ij}u_j).
\end{equation}
Furthermore, by shifting the time coordinate we may assume $0\in I$ and 
\begin{equation}
t_0=0.
\end{equation}
The partial differential equation in \re{5.38} can be expressed in the form
\begin{equation}\lae{5.41}
\begin{aligned}
&\msp[50]-a^{ij}u_{;ij}+\tilde v \bar h_{ij}\check u^i\check u^j=\\
&\msp[-40]e^\psi\{f(u,x)-\tilde\mu \bet_1(u-m_1)-\tilde\mu\bet_2(u-m_2)+n\al(e^{-\psi})_\al\tilde\nu^\al\}\\
&\msp[200]-v\bar H(u,x)
\end{aligned}
\end{equation}
or,  if we abbreviate the right-hand side of the  equation above by  
\begin{equation}\lae{5.42}
\tilde f=\tilde f(x,u,Du),
\end{equation}
in the form
\begin{equation}\lae{5.43}
-a^{ij}u_{;ij}+\tilde v \bar h_{ij}\check u^i\check u^j-\tilde f=0.
\end{equation}
Here,  the tensor $a^{ij}$ is defined by
\begin{equation}
a^{ij}=\tilde v\tilde  g^{ij}.
\end{equation}
Let $a^i$ be the vector field
\begin{equation}
a^i(p)=\tilde v \s^{ij}p_j,
\end{equation}
where
\begin{equation}
\tilde v=(1-\abs p^2)^{-1}
\end{equation}
and where we assume
\begin{equation}
\abs p^2=\s^{ij}p_ip_j<1,
\end{equation}
then
\begin{equation}
a^{ij}=\pde {a^i}{p_j}.
\end{equation}
To be absolutely precise, we should write
\begin{equation}
a^i=a^i(x,x^0,p)
\end{equation}
and correspondingly
\begin{equation}
a^{ij}=a^{ij}(x,x^0,p).
\end{equation}
In the equation \re{5.43} we consider
\begin{equation}
a^{ij}=a^{ij}(x,u,Du).
\end{equation}

Let $m_i'\in I$ be constants satisfying
\begin{equation}\lae{5.52}
m_1'<m_1-\e
\end{equation}
and
\begin{equation}\lae{5.53}
m_2'>m_2+\e,
\end{equation}
then we want to solve the variational inequality
\begin{equation}\lae{5.54}
\spd{-a^{ij}u_{;ij}+\tilde v \bar h_{ij}\check u^i\check u^j-\tilde f}{v-u}\ge 0\q\A\, v\in K,
\end{equation}
where the convex set $K$ is defined by
\begin{equation}
K=\{v\in H^{2,2}(B_R)\ii H^{1,2}_0(B_R):m_1'\le v\le m_2'\}
\end{equation}
and $u\in K$ is supposed to satisfy the additional requirement
\begin{equation}
\tilde v<\un.
\end{equation}
The duality in \re{5.54} is the real $L^2$ scalar product.
\br
Let us emphasize that for the definition of Sobolev  or H\"older spaces we equip the underlying domains, here, $B_R$, with a metric $\s_{ij}(x^0,x)$ with a constant $x^0$. In the present case we choose $x^0=0$. Note that the metrics $\s_{ij}(0,x)$ and $\s_{ij}(x^0,x)$ are all equivalent as long as $x^0$ stays in a compact subset of $I$.
\er

To find a solution of \re{5.54} we shall apply a version of the Leray-Schauder fixed point theorem. Consider in $C^{1,\al}_0(\bar B_R)$ the closed set
\begin{equation}
\msc C=\{w\in C^{1,\al}(\bar B_R):m_1'\le w\le m_2',\, \abs{Dw}^2\le 1-\theta,\, \abs w_{1,\al}\le c\},
\end{equation}
where
\begin{equation}
\abs{Dw}^2=\s^{ij}(w,x)w_iw_j,
\end{equation}
 and
\begin{equation}
C^{1,\al}_0(\bar B_R)=\{w\in C^{1,\al}(\bar B_R):\fv w{\pa B_R}=0\}.
\end{equation}
Here, $0<\al<1$, $\theta<1$ and $c$ are positive constants to be determined later.

$\msc C$ is closed with non-empty interior, since
\begin{equation}
0\in \inn{{\msc C}}.
\end{equation}
For $w\in \msc C$ consider the differential operator
\begin{equation}
Lu=-a^{ij}D_iD_ju+\tilde v\bar h_{ij}\check w^i\check u^j-\tilde f,
\end{equation}
where
\begin{equation}
a^{ij}=a^{ij}(x,w,Dw),
\end{equation}
\begin{equation}
\tilde v=(1-\abs{Dw}^2)^{-1},
\end{equation}
\begin{equation}
\s_{ij}=\s_{ij}(w,x)
\end{equation}
and \begin{equation}
\tilde f=\tilde f(x,w,Dw).
\end{equation}
The covariant differentiation is with respect to $\s_{ij}$. Define the map
\begin{equation}
T:\msc C\ra C^{1,\al}_0(\bar B_R)
\end{equation}
by the requirement that
\begin{equation}
u=Tw
\end{equation}
is a solution of the variational inequality
\begin{equation}\lae{5.68}
\spd{Lu-\tilde f}{v-u}\ge 0\qq\A\, v\in K.
\end{equation}
It is well-known that this variational inequality has a unique solution satisfying
\begin{equation}\lae{5.69}
Lu\in L^\un(B_R).
\end{equation}
Hence, we have
\begin{equation}\lae{5.70}
u\in H^{2,p}(B_R)\qq\A\, 1<p<\un,
\end{equation}
because of the Calderon-Zygmund inequalities. The existence of a solution of the variational inequality \re{5.68} and the relation \re{5.69} can be proved simultaneously with the help of the penalization method, i.e., by looking at Dirichlet problems similar to \re{5.14}, where in case of linear uniformly elliptic equations the Schauder theory guarantees the solvability of Dirichlet problems assuming appropriate smoothness of the data. A solution of the variational inequality can be achieved by letting the parameter $\e$, which enters in the definition of the penalization functions, tend to zero.

If we choose the exponent $p$ in \re{5.70} large enough such that
\begin{equation}
\bet=1-\frac np>\al,
\end{equation}
then the embedding
\begin{equation}
H^{2,p}(B_R)\hra C^{1,\al}(\bar B_r)
\end{equation}
is compact, hence $T$ is compact.

Now, let $u\in \msc C$ be an arbitrary quasi fixed point, i.e., there exists $\tilde\lam>1$ such that
\begin{equation}\lae{5.73}
Tu=\tilde\lam u.
\end{equation}
If we can show that
\begin{equation}
u\in \inn{\msc C},
\end{equation}
then $T$ will have a fixed point, \cf \cite[Theorem 4.43]{lloyd:book}, and it will be a solution of the variational inequality \re{5.54}. 

Thus, let $u$ be the quasi fixed point in \re{5.73}. In the open set
\begin{equation}
\Om=\{x\in B_r:m_1'<u(x)<m_2'\}
\end{equation}
it will solve the equation
\begin{equation}\lae{5.76}
-a^{ij}D_iD_ju+\tilde v \bar h_{ij}\check u^i\check u^j=\tilde\lam^{-1}\tilde f(x,u,Du).
\end{equation}
Adding on both sides
\begin{equation}
\tilde v^{-1}\bar H(x,u)
\end{equation}
we deduce, in view of \re{5.36}, that
\begin{equation}\lae{5.78} 
\tilde H =\tilde\lam^{-1}\tilde f (x,u,Du)+\tilde v^{-1}\bar H\qq\tup{in}\, \Om,
\end{equation}
where $\tilde H$ is the mean curvature of $M=\graph u$ in $\tilde N$. We  note that $\Om$  reaches up to the boundary because
\begin{equation}
m_1'<0<m_2',
\end{equation}
in view of \re{5.21} and the definition of the $m_i'$ in \re{5.52}and \re{5.53}, and that
\begin{equation}
\{x\in B_R: Du\not= 0\}\su\Om.
\end{equation}
We also emphasize that $u$, which is already of class $H^{2,p}(B_R)$, \cf \re{5.70}, satisfies
\begin{equation}
u\in C^{3,\al}(\Om\uu \pa B_R),
\end{equation}
because of the Schauder estimates, provided $\pa B_R$ is of class $C^{3,\al}$ which we assume.

By definition $u$ is bounded by the constants $m_i'$ and hence we can try to apply the a priori estimates in \frt{4.1} in the present situation. We only have to check that the additional terms on the right-hand side of \re{5.78} do not alter the structure of the inequality \fre{4.20} and  especially that the values of $\tilde\lam$ and $\bet_i'$ do not enter into the estimates. Recall that $\tilde f$ is defined to be the right-hand side of the equation \re{5.41}.

First, we observe that when we estimate $\tilde v$, or better,
\begin{equation}
w=\tilde v e^{\mu e^{\lam (u+c)}}
\end{equation}
we assume 
\begin{equation}
\tilde v >2,
\end{equation}
i.e., we consider points in $\Om$ for the interior estimates or points in $\pa B_R$ for the boundary estimates but never points where $Du=0$. 

In case of the boundary estimates we already switched to the conformal embedding so that the additional terms coming from that embedding are already taken care of, and from the original mean curvature $H$ only the supremum norm of $H$ will enter, hence only $\tilde\mu$ will enter into the estimates which is fine.

Thus, we only have to consider the interior estimates, where now the ambient metric is $\tilde g_{\al\bet}$. Let us look at the right hand side of \fre{3.20}, where now
\begin{equation}
(\h_\al)=(-1,0,\ldots, 0).
\end{equation}
The crucial term is the last one involving $f$ which is the prescribed mean curvature. In the present situation $f$ has to be replaced by
\begin{equation}\lae{5.85}
\tilde H=\tlam^{-1}\tilde f+\tilde v ^{-1}\bar H.
\end{equation}
The crucial term on the right hand-side of \re{3.20} has the form
\begin{equation}
-\tilde H_\bet x^\bet_i x^\al_k\tilde g^{ik}\h_\al=\tilde H_iu^i,
\end{equation}
where
\begin{equation}
u^i=\tilde g^{ij} u_j.
\end{equation}
From \re{5.85} we conclude that we only have to worry about
\begin{equation}
\tlam^{-1}\tilde f_iu^i
\end{equation}
and hence, in view of  \re{5.41} and  \re{5.42}, only about                                                       
\begin{equation}
-\tlam^{-1}\tilde\mu (\bet_1'+\bet_2')\norm{Du}^2\le 0
\end{equation}
which is non-positive and can be ignored, and in addition about the term
\begin{equation}
-n\tlam^{-1}\psi_{\al\bet}\tilde\nu^\al x^\bet_iu^i-n\tlam^{-1}\psi_\al x^\al_k\tilde h^k_iu^i,
\end{equation}
which can be estimated from above by
\begin{equation}
c\tilde v^3+c\norm{\tilde A}^2\norm{Du}^2\le c\tilde v^3+\e\norm{\tilde A}^2\tilde v+ c\e^{-1}\tilde v^3,
\end{equation}
where the $c$ on the right-hand side of the inequality is slightly larger than the $c$ on the left-hand side, and hence satisfies the structural requirements of the right-hand side of \fre{4.20}. Here, we also used
\begin{equation}
\tlam>1.
\end{equation}
Let us formulate the final conclusion as a lemma:
\bl\lal{5.3}
Let $u$ be a quasi solution of  $T$, then
\begin{equation}
\s^{ij}u_iu_j\le 1-\theta_0<1,
\end{equation}
where $\theta_0$ depends on the quantities mentioned in \frt{4.1} and furthermore on $\tilde\mu$ but not on $\tlam$ or $\bet_i'$. We also emphasize that $\abs u$ is bounded by 
$\max(m_2',-m_1')$ which can be made arbitrarily close to $\max(m_2+\e,-m_1-\e)$.
\el

Next, let us show that the $C^{1,\al}$-norm of $u$ is a priori bounded.
\bl\lal{5.4}
Let $u$ be a quasi fixed point of $T$, then
\begin{equation}
\abs{u}_{1,\al}\le c_1,
\end{equation}
where $c_1$ depends on $\theta_0$,  $\abs{\tilde f}_0$, $R$, $\max(m_2,-m_1)$ and the $C^2$-norm of $\pa B_R$ but is independent of the constants defining $\msc C$.
\el
\bp
Since we know that now  $a^{ij}$ is uniformly positive definite, let us write
\begin{equation}
-a^{ij}u_{;ij}=-D_i(\tilde v \s^{ij}D_ju)=-D_i(a^i(x,u,Du)),
\end{equation}
\cf \re{5.39}, where the covariant derivative is with respect to $\s_{ij}(u,x)$. From \re{5.69} we then obtain
\begin{equation}\lae{5.96}
-D_i(a^i(x,u,Du))=f^0\qq\tup{in}\, B_R,
\end{equation}
where
\begin{equation}
f^0\in L^\un(B_R)
\end{equation}
and that
\begin{equation}
\fv u{\pa B_R}=0.
\end{equation}
We know that $u$ is Lipschitz continuous, $\pa B_R$ of class $C^2$ and the vector field $a^i(x,x^0,p)$ of class $C^1$ in its arguments, at least for allowed values of $(x^0,p)$ and 
\begin{equation}
a^{ij}=\pde{a^i}{p_j}
\end{equation}
is uniformly elliptic for the allowed values. Hence it is known that
\begin{equation}
u\in H^{2,p}(B_R)\qq\A\, 1<p<\un
\end{equation}
and 
\begin{equation}
\norm u_{2,p}\le c_2,
\end{equation}
where $c_2$ depends on $p$, $R$, and known constants, for details see the appendix in \frs{9}.

Especially we have
\begin{equation}
\abs u_{1,\al}\le c_1,
\end{equation}
where $c_1$ is independent of the constants defining $\msc C$ only depending on the quantities mentioned in \rl{5.4}.
\ep

It finally remains to prove that a quasi fixed point of $T$ does not touch the obstacles.
\bl
Let $u$ be a quasi fixed point of $T$, then $u$ satisfies the inequalities
\begin{equation}
m_1'<u<m_2',
\end{equation}
and hence, $u$ also satisfies the equation \re{5.76} in $B_R$.
\el
\bp
We shall only prove the first inequality, since the proof of the other inequality is similar.

Let $x\in B_R$ be a point where
\begin{equation}
u(x_1)=m_1',
\end{equation}
then $u$ attains its minimum in $x_1$ and we have 
\begin{equation}
Du(x_1)=0.
\end{equation}
Moreover,
\begin{equation}
u<m_2'\qq\tup{in}\, B_\rho(x_1)
\end{equation}
for small $\rho$, hence we infer that in $B_\rho(x_1)$
\begin{equation}\lae{5.107} 
\begin{aligned}
&-a^{ij}D_iD_ju+\tilde v \bar h_{ij}\check u^i\check u^j\ge \tilde\lam^{-1}\tilde f(x,u,Du)\\
&= \tlam^{-1}e^\psi\{f(u,x)-\tilde\mu \bet_1(u-m_1)-n\psi_\al\tilde\nu^\al -v \bar H(u,x)\}.
\end{aligned}
\end{equation}
In $x_1$ we have
\begin{equation}
u-m_1=m_1'-m_1<-\e,
\end{equation}
hence, this inequality will also be valid in $B_\rho(x_1)$ if $\rho$ is small and we further deduce
\begin{equation}
\bet_1(u-m_1)=-1
\end{equation}
in $B_\rho(x_1)$. Thus, we conclude that in $x_1$
\begin{equation}\lae{5.110}
f(u,x)+\tilde\mu+n\dot\psi -\bar H(u,x)>0,
\end{equation}
if $\tilde\mu$ is large enough, which we stipulate. Since $u$ is of class $C^{1,\al}$ the right-hand side of  \re{5.107} will therefore be strictly positive if $\rho$ is small for the same $\tilde\mu$ as in \re{5.110}. But this leads to a contradiction when  we apply the weak Harnack inequality to $(u-m_1')$. Indeed, the left-hand side of \re{5.107} can be written in the form
\begin{equation}
-D_i(\tilde v\s^{ij}u_j)+\tilde v \bar h_{ij}\check u^i\check u^j
\end{equation}\and this uniformly elliptic divergence operator should be strictly positive in $B_\rho(x_1)$. But then $(u-m_1')$, which is non-negative,  cannot attain a minimum in $B_\rho(x_1)$ contrary to our assumption.
\ep

Thus, we have proved that $T$ has a fixed point which is actually a solution of the Dirichlet problem \fre{5.14}, and hence, also a solution of the original Dirichlet problem \fre{3.1}, in view of \rt{5.1}.

\section{Decay estimates for the solutions}

In this section we want to derive decay estimates for solutions of the Dirichlet problems in \fre{3.1}, where, now and for the rest of the paper, we assume for simplicity that
\begin{equation}
f(x^0,x)=\hat H(t)\qq\, t\in I,
\end{equation}
i.e., in case of the Dirichlet problem \re{3.1}
\begin{equation}
f=\hat H(t_0),
\end{equation}
for otherwise we would have to impose appropriate conditions on $f$. $\hat H(t)$ is the mean curvature of the slice
\begin{equation}
\{x^0=t\}
\end{equation}
in the Robertson-Walker spacetime $\hat N$.  The decay parameter is $r$, the radial distance function; $r$ is the geodesic distance from a fixed point of the underlying space of constant curvature, i.e., $\R[n]$ or $\Hh[n]$, where, to be absolutely precise, the metric of these spaces is multiplied by the scale factor $a^2(x^0)$, and hence, the corresponding distance $r$ also depends on $a(x^0)$ 
\begin{equation}\lae{6.1}
r=a(x^0)\rho,
\end{equation}
where $\rho$ is the geodesic distance corresponding to $a=1$.

For the decay estimates we use the radial function in \re{6.1} with
\begin{equation}
x^0=t_0,
\end{equation}
the boundary values of the corresponding Dirichlet problem. Note, that the balls
\begin{equation}
B_R\su\so
\end{equation}
are defined with respect to $\rho$ in order to have a common domain of definition in case we want to compare solutions with different data.

For simplicity we shall redefine the balls for the decay estimates by assuming that the radius $R$ is defined by the geodesic distance
\begin{equation}
r=a(t_0)\rho.
\end{equation}
We then look at the Dirichlet problem \re{3.1} in balls with $R>R_0$, where $R_0$ is so large such that the asymptotical estimates in \fras{2.3} are valid, especially the estimate in \re{2.52}, namely,
\begin{equation}\lae{6.5}
\dot{\bar H}\ge c_0>0
\end{equation}
uniformly in
\begin{equation}
\{(x^0,x):r(x)\ge R_0,\, m_1\le x^0\le m_2\},
\end{equation}
where the $m_i$, $i=1,2$, are bounds for the solutions $u$ of the Dirichlet problems \re{3.1} which are independent of $R$.

Then can prove:
\bl\lal{6.1}
Let $u$ be a solution of the Dirichlet problem \re{3.1}, then there exists $R_1>R_0$ such that $\abs{u-t_0}$ can be estimated by
\begin{equation}
\abs{u(x)-t_0}\le \frac cr\q\A\, x\in \bar B_R\ii\{r(x)\ge R_1\},
\end{equation}
where $R_1>R_0$ and the constants $R_1$ and $c$ depend on $c_0$ in \re{6.5} and other already known estimates for $\abs{Du}$ and $\abs u$ as well as the ambient capped region $\msc S(m_1,m_2)$.
\el
\bp
We shall construct appropriate barrier functions to estimate $u$ from above and from below. Let $\f=\f(r)$ be defined by 
\begin{equation}
\f(r)=\frac{\lam R_1}r+t_0
\end{equation}
be the upper barrier and
\begin{equation}
\tilde\f(r)=-\frac{\lam R_1}r+t_0
\end{equation}
be the lower barrier, where $\lam>0$ is a fixed constant satisfying
\begin{equation}
\lam>\sup \abs u +\abs{t_0}.
\end{equation}
The value of $R_1$ will be determined later.  Since the arguments for the upper \resp lower estimates are similar we shall only prove the upper estimate
\begin{equation}\lae{6.11}
u(x)\le \frac{\lam R_1}r+t_0\qq \A\, x\in B_R\ii\{r(x)\ge R_1\}.
\end{equation}
Let $\Om$ be the region
\begin{equation}
\Om=\{x\in B_R:R_1<r(x)<R\},
\end{equation}
then we to prove the estimate \re{6.11} in $\Om$. The boundary of $\Om$ consists of the spheres $\{r(x)=R_1\}$ and $\{r(x)=R\}$. On both spheres the inequality \re{6.11} is strict. Thus, let us argue by contradiction and assume that there exists $x_0\in\Om$ such that
\begin{equation}
u(x_0)-\f(x_0)=\sup_\Om(u-\f)>0.
\end{equation}
Then we know that in $x_0$
\begin{equation}
Du=D\f
\end{equation}
and
\begin{equation}
u_{,ij}\le \f_{,ij},
\end{equation}
where the comma indicates partial derivatives. In order to apply the maximum principle we again consider $M$ to be embedded in $(\tilde N,\tilde g_{\al\bet})$, then $u$ satisfies the elliptic equation
\begin{equation}\lae{6.16} 
\tilde H \tilde v=-\tilde v^2\tilde g^{ij}u_{;ij}+\tilde{\bar H}+\tilde v^2\tilde{\bar h}_{ij} \check u^i\check u^j,
\end{equation}
where $\tilde{\bar h}_{ij}$ is the second fundamental form of the slices
\begin{equation}
\{x^0=\const\}
\end{equation}
in $\tilde N$, \cf \fre{5.36}. The geometric quantities are related to the corresponding geometric quantities in $N$ by
\begin{equation}
e^\psi H=\tilde H+n\psi_\al\tilde\nu^\al
\end{equation}
and
\begin{equation}
e^\psi \bar H=\tilde{\bar H}-n \dot\psi.
\end{equation}
Furthermore, we know that
\begin{equation}
\abs{e^\psi-1}+\nnorm{D\psi}\le \frac cr\qq\A\, x\in \{r(x)\ge R_0\}
\end{equation}
and
\begin{equation}\lae{6.24}
\tilde v-1=\frac{\tilde v^2}{\tilde v+1}\abs{Du}^2
\end{equation}
Thus, we deduce
\begin{equation}\lae{6.21}
\hat H(t_0)\equiv f(t_0)=H\ge -\tilde v^2\tilde g^{ij}u_{;ij}+\bar H(u,x)-c \abs{Du}^2-\frac cr\qq\A\,x\in\Om,
\end{equation}
where $\hat H(t_0)$ is the mean curvature of the slice 
\begin{equation}
\{x^0=t_0\}
\end{equation}
in the Robertson-Walker spacetime $\hat N$ and the identities on the left-hand side of the inequality simply reflect the prescribed value of $H$.

Next, let us consider the covariant derivatives of $u$. We have
\begin{equation}\lae{6.23} 
\begin{aligned}
-u_{;ij}&=-u_{,ij}+\C^k_{ij}u_k\\
&\ge -\f_{,ij}+\hat\C^k_{ij}\f_k+(\C^k_{ij}-\hat\C^k_{ij})\f_k,
\end{aligned}
\end{equation}
where $\hat\C^k_{ij}$ are the Christoffel symbols of the metric $\hat \s_{ij}$ which is the induced metric of the slice
\begin{equation}
\{x^0=t_0\}
\end{equation}
in $\hat N$. In view of the boundedness of $\tilde v$ and the assumptions \re{2.47} and \fre{2.48} we infer that the tensor
\begin{equation}
\C^k_{ij}-\hat\C^k_{ij}
\end{equation}
is bounded relative to the metric $\s_{ij}$.

Moreover,
\begin{equation}
\f_i=-\frac{\lam R_1}{r^2}r_i,
\end{equation}
\begin{equation}
\f_{,ij}=2\frac{\lam R_1}{r^3}r_ir_j-\frac{\lam R_1}{r^2}r_{,ij}
\end{equation}
and the tensor
\begin{equation}
r_{,ij}-\hat\C^k_{ij}r_k
\end{equation}
is uniformly bounded with respect to $\hat \s_{ij}$ and hence also with respect to $\s_{ij}$. Thus, we conclude from \re{6.21} and \re{6.23}
\begin{equation}
\hat H(t_0)\ge \bar H(u,x)-\frac cr
\end{equation}
in $x=x_0$, where $c$ is independent of $R_1$. Furthermore, in view of \fre{2.47},
\begin{equation}
\hat H(t_0)\le \bar H(t_0,x)+\frac cr\qq\A\, x\in\Om
\end{equation}
and therefore
\begin{equation}
\bar H(t_0,x_0)\ge \bar H(u,x)-\frac cr,
\end{equation}
or equivalently,
\begin{equation}
0\ge \bar H(u,x)-\bar H(t_0,x_0)-\frac cr,
\end{equation}
from which we conclude
\begin{equation}
\begin{aligned}
0&\ge \dot{\bar H}(\tau,x_0)(u-t_0)-\frac cr\\
&\ge c_0\frac{\lam R_1}r-\frac cr>0,
\end{aligned}
\end{equation}
in view of \fre{2.52}, provided
\begin{equation}
R_1>\frac c{\lam c_0},
\end{equation}
a contradiction.
\ep

The decay estimate for $\abs{u-t_0}$ allows us to prove a similar decay estimate for $\abs{Du}$ and $\abs{D^2u}$, and also for higher derivatives, with the help of the Schauder estimates. But first, let us derive local Schauder estimates.

\bl\lal{6.2}
Let $u$ be a solution of the Dirichlet problem \re{3.1} of class $C^{3,\al}$, then $u$ satisfies the elliptic differential equation
\begin{equation}\lae{6.38}
-a^{ij}u_{;ij}+b^iu_i+c(x) (u-t_0)=\tilde f,
\end{equation}
where the covariant differentiation is with respect to $\s_{ij}(u,x)$ and
\begin{equation}
a^{ij}=\tilde v^2\tilde g^{ij},
\end{equation}
\begin{equation}
\begin{aligned}
b^i&=\tilde v^2\dot\psi \s^{ij}u_j-n \tilde v^2\psi_j\s^{ji}-(n-1)\dot \psi\tilde v^2 \s^{ij}u_j\\
&\hp{=}+\tilde v^2 e^{\psi}\bar h_{jk}\check u^j\s^{ik}-e^\psi \hat H(t_0)\frac{\tilde v^2}{\tilde v+1}\s^{ji}u_j,
\end{aligned}
\end{equation}
\begin{equation}\lae{6.41}
c=e^\psi\int_0^1\dot{\bar H}(tu+(1-t)t_0,x)dt
\end{equation}
and
\begin{equation}
\tilde f=e^\psi(\hat H(t_0)-\bar H(t_0,x)).
\end{equation}
Then, for any $x_0\in B_R$ and any $\rho_0>0$ such that
\begin{equation}
\bar B_{2\rho_0}(x_0)\su B_R
\end{equation}
the estimate
\begin{equation}\lae{6.44.1} 
\abs{u-t_0}_{2,\al,B_{\rho_0}(x_0)}\le c (\abs{\tilde f}_{0,\al,\B_{2\rho_0}(x_0)}+\abs{u-t_0}_{0,B_{2\rho_0}(x_0)})
\end{equation}
is valid, where $c$ depends on $\rho_0$ and known quantities but not on $R$.

Similarly, for any $x_0\in\pa B_R$ and any $\rho_0>0$ we have
\begin{equation}\lae{6.45.1}
\abs{u-t_0}_{2,\al, \Om(x_0,\rho_0)}\le c(\abs{\tilde f}_{0,\al,\Om(x_0,2\rho_0)}+\abs{u-t_0}_{0,\Om(x_0,2\rho_0)}),
\end{equation}
where $c$ depends on $\rho_0$ and known quantities but not on $R$. Here,
\begin{equation}
\Om(x_0,\rho_0)=B_R\ii B_{\rho_0}(x_0).
\end{equation}
If the ambient space $N$ is of class $C^{m,\al}$, $m\ge 3$, $0<\al<1$, then
\begin{equation}\lae{6.47.1}
\abs{u-t_0}_{m,\al,B_{\rho_0}(x_0)}\le c (\abs{\tilde f}_{m-2,\al,\B_{2\rho_0}(x_0)}+\abs{u-t_0}_{0,B_{2\rho_0}(x_0)})
\end{equation}
where $c$ depends on $\rho_0$, $m$, and known quantities but not on $R$, and near the boundary
\begin{equation}\lae{6.48.1}
\abs{u-t_0}_{m,\al,\Om(x_0,\rho_0)}\le c (\abs{\tilde f}_{m-2,\al,\Om(x_0,2\rho_0)}+\abs{u-t_0}_{0,\Om(x_0,2\rho_0)})
\end{equation}
for any $x_0\in \pa B_R$ and $\rho_0>0$, where $c$ depends on $\rho_0$,  $m$  and known quantities but not on $R$.
\el

\bp
\cq{\re{6.38}}\q Follows immediately from \re{6.16} and \re{6.24}.

\cvm
\cq{\re{6.44.1}}\q Since the coefficients in the operator in \re{6.38} are uniformly of class $C^{0,\al}$, \cf \frr{9.1}, the estimate follows immediately from the interior Schauder estimates, \cf \cite[Corollary 6.3]{gilbarg}.

\cvm
\cq{\re{6.45.1}}\q This estimate is due to the Schauder estimates near a $C^{2,\al}$ boundary having in mind that $(u-t_0)$ has zero boundary values, \cf \cite[Corollary 6.7]{gilbarg}.

\cvm
\cq{\re{6.47.1}}\q Let us express the covariant derivatives in equation \re{6.38} as partial derivatives and let us differentiate this equation with respect to $x^k$, where $1\le k\le n$ is fixed. The resulting equation is an  elliptic equation for 
\begin{equation}
w=D_k(u-t_0).
\end{equation}
Applying then the estimate \re{6.44.1} to $w$, while replacing $\rho_0$ by $\frac{\rho_0}2$, we deduce
\begin{equation}
\begin{aligned}
\abs w_{2,\al,B_{\frac{\rho_0}2}(x_0)}&\le c(\abs {\tilde f}_{1,\al,B_{\rho_0}(x_0)}+\abs{u-t_0}_{2,\al,B_{\rho_0}(x_0)})\\
&\le c(\abs {\tilde f}_{1,\al,B_{2\rho_0}(x_0)}+\abs{u-t_0}_{0,B_{2\rho_0}(x_0)}),
\end{aligned}
\end{equation}
in view of the previous result for $(u-t_0)$. Here the quotient of the radii is $4$ instead of $2$, but it is obvious that any number larger than $1$ would suffice.

Since $k$ is arbitrary, we have the estimate for $m=3$. The estimate for larger $m$ is then proved inductively having in mind that $\rho_0$ is flexible.

\cvm
\cq{\re{6.48.1}}\q For the $C^{m,\al}$-estimates near the boundary we have to straighten the boundary either explicitly, or implicitly, by using geodesic polar coordinates $(r,\xi^a)$, $1\le a\le n-1$, with respect to the induced metric of the slice
\begin{equation}
\{x^0=t_0\}
\end{equation}
in $\hat N$. This coordinate system is also an allowed coordinate system if the slice is embedded  in $N$ though not a  geodesic polar coordinate system. But the coordinates $(R,\xi^a)$ describe the boundary $\pa B_R$ which is all we need. The partial derivatives with respect $\xi^a$ are then tangential derivatives and the function
\begin{equation}
w=D_a(u-t_0)
\end{equation}
has zero boundary values. Similarly as in \re{6.45.1} we obtain, for any $x_0\in B_R$ and any $0<\rho_0<\frac R4$,
\begin{equation}\lae{6.53.1}
\begin{aligned}
\abs w_{2,\al,\Om(x_0,{\frac{\rho_0}2})}&\le c(\abs {\tilde f}_{1,\al,\Om(x_0,\rho_0)}+\abs{u-t_0}_{2,\al,\Om(x_0,\rho_0)})\\
&\le c(\abs {\tilde f}_{1,\al,\Om(x_0,2\rho_0)}+\abs{u-t_0}_{0,\Om(x_0,2\rho_0)}),
\end{aligned}
\end{equation}
where $c$ is independent of $R$.

Let us denote the radial coordinate by $\xi^n$, then, in \re{6.53.1}, we have estimated
\begin{equation}\lae{6.54.1}
\sum_{i+j+k<3n}\abs{D_iD_jD_k (u-t_0)}_{0,\al,\Om(x_0,{\frac{\rho_0}2})}.
\end{equation}
The remaining derivative
\begin{equation}
D_nD_nD_nu
\end{equation}
can be estimated by looking at the uniformly elliptic equation for 
\begin{equation}
w=D_nu
\end{equation}
and applying the estimate for \re{6.54.1} as well as the uniform positivity of $a^{nn}$.

The $C^{m,\al}$-estimates near the boundary for $m>3$ can then be proved inductively.
\ep

Combining the results of \rl{6.1}, \rl{6.2}  and the decay assumptions for $(e^\psi-1)$ and its derivatives \resp  $(\hat H(t_0)-\bar H(t_0,x))$ and its  spatial derivatives, \cf \fras{2.3},  we can  state:
\bt\lat{6.3}
Let $N$ be of class $C^{m,\al}$, $m\ge 3$, $0<\al<1$, then for any ball
\begin{equation}
\bar B_{2\rho_0}(x_0)\su B_R\sminus B_{R_1}
\end{equation}
\resp boundary neighbourhood
\begin{equation}
\Om(x_0,2\rho_0)\su \bar B_R\sminus B_{R_1}
\end{equation}
the decay estimates
\begin{equation}
\abs{u-t_0}_{m,\al,B_{\rho_0}(x_0)}\le \frac c{r(x_0)}
\end{equation}
\resp
\begin{equation}
\abs{u-t_0}_{m,\al,\Om(x_0,\rho_0)}\le \frac c{r(x_0)}
\end{equation}
are valid, where $c$ depends on $m$, $\rho_0$ and known quantities but not on $R$. Especially we obtain
\begin{equation}\lae{6.61} 
\sum_{0\le \abs\ga\le m}\abs{D^\ga(u-t_0)(x)}\le \frac c{r(x)}\qq\A\, x\in \bar B_R\sminus B_{R_1},
\end{equation}
where $c$ depends on $m$ but not on $x$ or $R$.
\et
\bp
Obvious.
\ep

\bc\lac{6.4}
Let $u=u_R$ be a solution of the Dirichlet problem \re{3.1}, then there exists $c_0>0$ such that
\begin{equation}\lae{6.62}
(\abs A^2+\bar R_{\al\bet}\nu^\al\nu^\bet)\ge c_0>0\qq\A\,x\in \{r(x)>R_1\},
\end{equation}
if $R_1$ is large enough.
\ec
\bp
This is due to the estimate \re{6.61} which implies that the hypersurface
\begin{equation}
M=\graph u
\end{equation}
converges to the slice
\begin{equation}
M(t_0)=\{x^0=t_0\}
\end{equation}
in the class $C^2$ if $r(x)$ is large, i.e., the left-hand side of \re{6.62} uniformly converges to the corresponding quantity of $M(t_0)$ which satisfies this inequality because of \re{2.52}, \re{2.49} and the relation
\begin{equation}
\dot H=-\D e^\psi+(\abs A^2+\bar R_{\al\bet}\nu^\al\nu^\bet)e^\psi
\end{equation}
which is valid for the slices
\begin{equation}
\{x^0=\const\},
\end{equation}
\cf \cite[equ. (2.3.37), p. 96]{cg:cp}.
\ep

\section{Existence of a CMC foliation}\las{7}

Let $u_R$ be a solution of the Dirichlet problems \fre{3.1} with boundary values $t_0$. In view of the estimates in the previous section we infer that, by letting $R$ tend to infinity, a subsequence converges in
\begin{equation}
C^{m,\al}_{\tup{loc}}(\so)\ii C^m_{\tup{loc}}(\bar \so),
\end{equation}
where the subscript ``loc'' is necessary since the $u_R$ are only defined in $B_R$, to a function
\begin{equation}\lae{7.2}
u\in C^{m,\al}(\so)\ii C^m(\bar \so)
\end{equation}
such that $M=\graph u$ is a spacelike hypersurface of constant mean curvature
\begin{equation}\lae{7.3}
\fv FM=f(t_0)=\hat H(t_0).
\end{equation}
Here, the function space
\begin{equation}
C^m(\bar\so)\su C^m(\so)
\end{equation}
is defined by the additional requirements
\begin{equation}
t=\lim_{r(x)\ra \un} u(x)\in I\qq\A\, u\in C^2(\bar\so)
\end{equation}
and that
\begin{equation}\lae{7.6}
\abs{u(x)-t}\le \frac c{r(x)}\qq\A\, x\in \{r(x)\ge R_1\}
\end{equation}
as well as
\begin{equation}
\sum_{0\le \abs\ga\le m}\abs{D^\ga(u-t_0)(x)}\le \frac c{r(x)}\qq\A\, x\in \{r(x)\ge R_1\}
\end{equation}
for suitable large constants $c$ and $R_1$.

Let us write
\begin{equation}\lae{7.8}
u=u(t_0,x)
\end{equation}
for the function in \re{7.2}, where
\begin{equation}
t_0=\lim_{r(x)\ra\un}u(x).
\end{equation}
Since the mean curvature function
\begin{equation}
\tau=\hat H(t)\qq t\in I
\end{equation}
is invertible, recall that
\begin{equation}
\dot{\hat H}>0,
\end{equation}
we could also express $u$ in the form
\begin{equation}\lae{7.12}
u=u(\tau_0,x),
\end{equation}
though, at the moment, the convention \re{7.8} is more suitable. We shall also speak of \tit{the} solution, since we shall later prove that the sequence $u_R$ actually converges and not only subsequences, because the solutions of equation \re{7.3} satisfying the estimate \re{7.6} are uniquely determined.

Moreover, in view of the Comparison \frt{3.3} we conclude that
\begin{equation}\lae{7.13}
t_1<t_2\ra u(t_1,x)<u(t_2,x)\qq\A\, (t_i, x)\in I\times \so.
\end{equation}
The strict inequality on the right-hand side is again due to the weak Harnack inequality, \cf the proof of \rl{7.2}, where this argument is applied in a more detailed fashion.

The functions $u(t,\cdot)$  also satisfy:
\bl\lal{7.1}
Let 
\begin{equation}
M_t=\graph u(t,\cdot)\qq t\in I=(t_-,t_+),
\end{equation}
then
\begin{equation}\lae{7.15} 
\lim_{t\ra t_+}u(t,x)=t_+\qq\A\, x\in \so
\end{equation}
and
\begin{equation}
\lim_{t\ra t_-}u(t,x)=t_-\qq\A\, x\in \so.
\end{equation}
\el
\bp
We shall only prove the first relation. Let $M_k=\graph \psi_k$ be an arbitrary barrier such that
\begin{equation}
H_k=\fv H{M_k}
\end{equation}
and
\begin{equation}
t_k=\sup \psi_k,
\end{equation}
then, if $u_R(t,\cdot)$ is a solution of the Dirichlet problem in $B_R$ with boundary value $t$ and mean curvature $\hat H(t)$, we deduce
\begin{equation}\lae{7.19}
u_R(t,x) >\psi_k(x)\qq\A\, x\in B_R
\end{equation}
provided
\begin{equation}\lae{7.20}
t>t_k\q\wed\q \hat H(t)>H_k,
\end{equation}
in view of \frt{3.3}. The requirements \re{7.20} can be easily satisfied because
\begin{equation}
\lim_{t\ra t_+}\hat H(t)=H_+.
\end{equation}
The estimate \re{7.19} is then also valid for $u(t,\cdot)$, provided \re{7.20} is true, from which  the relation \re{7.15} can be easily deduced in view of the properties of the barriers.
\ep

Let us now show that the spacelike hypersurfaces over $\so$
\begin{equation}
M=\graph u
\end{equation}
satisfying
\begin{equation}\lae{7.23}
\fv HM=c_0=\const
\end{equation}
and
\begin{equation}\lae{7.24}
\abs{u(x)-t_0}\le \frac c{r(x)}\qq\A\, x\in \{r(x)>R_1\}
\end{equation}
are uniquely determined.
\bl\lal{7.2}
Let $M=\graph u$ be a spacelike hypersurface of class $C^{3,\al}(\so)$ satisfying \re{7.23} and \re{7.24}, then $M$ is unique.
\el
\bp
We argue by contradiction by assuming there exists two  hypersurfaces
\begin{equation}
M_i=\graph u_i,\qq i=1,2,
\end{equation}
satisfying \re{7.23} as well as \re{7.24}. Let us suppose that
\begin{equation}
0<d_0=d(M_2,M_2),
\end{equation}
then there exists $p_k\in M_2$ and $q_k\in M_1$ such that
\begin{equation}
d(p_k,q_k)\ra d_0.
\end{equation}
The points can be expressed in the form
\begin{equation}
p_k=(u_2(x_k),x_k)
\end{equation}
and similarly
\begin{equation}
q_k=(u_1(y_k),y_k).
\end{equation}
If $r(x_k)$ remains uniformly bounded, then $r(y_k)$ is also uniformly bounded because
\begin{equation}
q_k\in I^+(p_k),
\end{equation}
and vice versa. On the other hand, if $r(x_k)$ and $r(y_k)$ both tend to infinity, then for any $\e>0$ there exists $k_0$ such that
\begin{equation}
t_0-\e<u_1(y_k)<t_0+\e\qq\A\, k\ge k_0
\end{equation}
and
\begin{equation}
t_0-\e<u_2(x_k)<t_0+\e\qq\A\, k\ge k_0.
\end{equation}
Let $\ga_k$ be the future directed curve connecting $p_k$ and $q_k$, then $\ga_k$ can be extended to a future directed curve $\tilde \ga_k$ connecting the slices
\begin{equation}
M(t_0-\e)=\{x^0=t_0-\e\}
\end{equation}
and 
\begin{equation}
M(t_0+\e)=\{x^0=t_0+\e\}
\end{equation}
such that the lengths of the extended curves satisfy
\begin{equation}
L(\tilde\ga_k)>L(\ga_k)\ge \frac {d_0}2\qq\A\, k\ge k_0,
\end{equation}
if $k_0$ is large enough. But this is a contradiction since then
\begin{equation}
\frac {d_0}2\le d(M(t_0-\e),M(t_0+\e))\equiv d_\e
\end{equation}
and
\begin{equation}
\lim_{\e\ra0}d_\e=0.
\end{equation}
Hence, we conclude that there exist $p_i$, $i=1,2$,
\begin{equation}
p_i\in M_i
\end{equation}
such that
\begin{equation}\lae{7.39}
d_0=d(p_2,p_1)=d(M_2,M_1),
\end{equation}
where $d$ is the Lorentzian distance function. Let $\f$ be a maximal geodesic from
$M_2$ to $M_1$  realizing this distance with endpoints $p_2$ and $p_1$, and
parametrized by arc length.

Denote by $\bar d$ the Lorentzian distance function to $M_2$, i.e., for $p\in
I^+(M_2)$
\begin{equation}
\bar d(p)=\sup_{q\in M_2}d(q,p).
\end{equation}
Since $\f$ is maximal, $\C=\set{\f(t)}{0\le t<d_0}$ contains no focal points of
$M_2$,
\cf \cite[Theorem 34, p. 285]{bn}, hence there exists an open neighbourhood $\mc
V=\mc V(\C)$ such that $\bar d$ is smooth in $\mc V$, \cf \cite[Proposition
30]{bn}, because $\bar d$ is a component of the inverse of the normal exponential
map of $M_2$.

%\cvm
Now, $M_2$ is the level set $\{\bar d=0\}$, and the level sets
\begin{equation}
M(t)=\set{p\in \mc V}{\bar d(p)=t}
\end{equation}
are  $C^3$ hypersurfaces; $x^0=\bar d$ is a time function in $\mc V$ and 
generates a normal Gaussian coordinate system, since $\spd{D\bar d}{D\bar d}=-1$.
Thus, the mean curvature $\bar H(t)$ of $M(t)$ satisfies the equation
\begin{equation}
\dot {\bar H}=\abs{\bar A}^2+\bar R_{\al\bet}\nu^\al\nu^\bet,
\end{equation}
\cf \cite[equ. (2.3.27, p. 96]{cg:cp}, and therefore we have
\begin{equation}\lae{2.15b}
\dot{\bar H}\ge 0,
\end{equation}
in view of the assumption \re{2.53}.

%\cvm
Next, consider a local tubular neighbourhood $\mc U$ of $M_1$ near $p_1$---to simply the phrasing---with corresponding
normal Gaussian coordinates $(x^\al)$. The (local) level sets
\begin{equation}
\tilde M(s)=\{x^0=s\},\qq-\e<s<0,
\end{equation}
lie in the past of $M_1=\tilde M(0)$ and are of class $C^{3,\al}$ for small $\e$.

Since the geodesic $\f$ is normal to $M_1$, it is also normal to $\tilde M(s)$ and
the length of the geodesic segment of $\f$ from $\tilde M(s)$ to $M_1$ is exactly
$-s$, i.e., equal to the distance from $\tilde M(s)$ to $M_1$, hence we deduce
\begin{equation}\lae{7.45}
d(M_2,\tilde M(s))=d_0+s,
\end{equation}
i.e., $\set{\f(t)}{0\le t\le d_0+s}$ is also a maximal geodesic from $M_2$ to $\tilde
M(s)$, and we conclude further that, for fixed $s$, the hypersurface $\tilde
M(s)\ii\mc V$ is contained in the past of $M(d_0+s)$ and touches $M(d_0+s)$ in
$p_s=\f(d_0+s)$. 

Let $p_2$ have the coordinates $(0,x_0)$ in the normal Gaussian coordinate system defined in the tubular neighbourhood $\mc V$, then there exists a small ball around $x_0$
\begin{equation}
B_\de(x_0)\su M_2
\end{equation}
such that the cylinder
\begin{equation}
Q=[0,d_0+s]\times B_\de(x_0)\su \mc V
\end{equation}
and the hypersurface
\begin{equation}
\tilde M(s)\ii  Q
\end{equation}
can be written as a graph over $B_\de(x_0)$, where the time coordinate is $x^0=\bar d$. Note that
\begin{equation}
\hat M(s)=\tilde M(s)\ii Q=\{(x^0,x):x^0=u(x),\; x\in B_\de(x_0)\}
\end{equation}
lies in the past of the slice $M(d_0+s)$ in view of \re{7.45}, hence if we write
\begin{equation}
\hat M(s)=\graph u
\end{equation}
in $B_\de(x_0)$, then
\begin{equation}
u\le u(x_0)=d_0+s.
\end{equation}
The hypersurface $\hat M(s)$ therefore touches the hypersurface $M(d_0+s)$ from below. Since
\begin{equation}
\bar H(d_0+s)=\fv H{M(d_0+s)}\ge c_0
\end{equation}
and
\begin{equation}
\fv H{\hat M(s)}\le c_0
\end{equation}
we deduce
\begin{equation}
\bar H(d_0+s)-\fv H{\hat M(s)}\ge 0
\end{equation}
in $B_\de(x_0)$. The left-hand side can be written as a uniformly elliptic equation for the non-negative function
\begin{equation}
\tilde u=d_0+s-u\ge 0
\end{equation}
such that
\begin{equation}
-D_i(a^{ij}D_j\tilde u)+b^i\tilde u_i+c\tilde u\ge 0
\end{equation}
in $\B_\de(x_0)$, where we used \fre{6.16}, \fre{5.39}, \fre{6.41} and the identity
\begin{equation}
\tilde v-1 = \frac{\tilde v^2}{\tilde v+1}\abs {Du}^2.
\end{equation}
The weak Harnack inequality then implies
\begin{equation}
u(x)=d_0+s\qq\A\, x\in B_\de(x_0).
\end{equation}
But then the tubular neighbourhood contains the full cylinder
\begin{equation}
[0,d_0]\times B_\de(x_0)
\end{equation}
and for every point
\begin{equation}
p=(0,x)\in M_2,\qq\, x\in B_\de(x_0),
\end{equation}
there exists a corresponding point $q\in M_1$ such that
\begin{equation}
d_0=d(M_2,M_1)=d(p,q).
\end{equation}
Now we can derive a contradiction. Let $\Lam$ be defined by
\begin{equation}
\Lam=\{p\in M_2:\E \, q\in M_1\,\tup{with}\, d(p,q)=d_0\}.
\end{equation}
Obviously, $\Lam\not=\eS$ and closed and we have just proved that $\Lam$ is also open, hence
\begin{equation}
\Lam=M_2
\end{equation}
and any future directed orthogonal geodesic emanating from $M_2$ will meet $M_1$ and realize the distance $d_0$, which  certainly contradicts the estimate \re{7.24}.
\ep

\bc\lac{7.3}
The functions $u(t,\cdot)$ the graphs of which have constant mean curvature $\hat H(t)$ and which tend to $t$ at spatial infinity are continuous in $t$.
\ec
\bp
Obvious.
\ep

Thus, there exists a family
\begin{equation}
M_t=\graph u(t,\cdot),\qq t\in I,
\end{equation}
of constant mean curvature hypersurfaces which are monotonically ordered by \re{7.13}. We shall now prove that they form a foliation of $N$.
\bt\lat{7.4}
The hypersurfaces $M_t$, $t\in I$, provide a foliation of $N$.
\et
\bp
We argue by contradiction and assume that there exists a point $p_0=(t_0,x_0)$ such that
\begin{equation}
p_0\notin \uuu_{t\in I}M_t.
\end{equation}
Let $\Lam_\pm$ be defined by
\begin{equation}
\Lam_+=\{t\in I:u(t,x_0)>t_0\}
\end{equation}
and
\begin{equation}
\Lam_-=\{t\in I:u(t,x_0)<t_0).
\end{equation}
Both sets are nonempty because of \rl{7.1}. 

Next, let $u_2$ be defined by
\begin{equation}
u_2=\inf_{t\in \Lam_+}u(t,\cdot)
\end{equation}
and
\begin{equation}
u_1=\sup_{t\in \Lam_-}u(t,\cdot),
\end{equation}
then
\begin{equation}
\lim_{r(x)\ra\un}u_2(x)=t_2=\inf_{t\in\Lam_+}t
\end{equation}
and
\begin{equation}
\lim_{r(x)\ra\un}u_1(x)=t_1=\sup_{t\in\Lam_-}t
\end{equation}
in view of \re{7.24} and the respective mean curvatures of 
\begin{equation}
M_i=\graph u_i,\qq\, i=1,2,
\end{equation}
are $\hat H(t_i)$. Moreover, because of the monotonicity
\begin{equation}
t_1\le t_2,
\end{equation}
\begin{equation}
u_1\le u_2
\end{equation}
and
\begin{equation}
u_1(x_0)\le t_0\le u_2(x_0).
\end{equation}
Because of \rl{7.2} we know
\begin{equation}
u_i=u(t_i,\cdot).
\end{equation}
Thus, the assumption
\begin{equation}
t_1<t_2
\end{equation}
would lead to a contradiction by simply picking any $\bar t$
\begin{equation}
t_1<\bar t<t_2.
\end{equation}
The existence of $\graph u(\bar t,\cdot)$ would then contradict the definition of $u_1$ or $u_2$. Hence, we conclude
\begin{equation}
t_1=t_2
\end{equation}
and consequently
\begin{equation}
u_1=u_2
\end{equation}
which implies
\begin{equation}
u_1(t_0)=t_0,
\end{equation}
i.e., the theorem is proved.
\ep

\section{The constant mean curvature $\tau$ of the foliation by spacelike hypersurfaces is a smooth global  time function}

After having proved the existence of a unique foliation by spacelike CMC  hypersurfaces $M_t$, $t\in I$, we want to show that the mean curvature
\begin{equation}\lae{8.1}
\tau=\hat H(t)
\end{equation}
of the foliation hypersurfaces is a smooth global time function, provided the  spacetime $N$ is smooth. First, $\tau$ is certainly a continuous time function, since the level hypersurfaces
\begin{equation}
\{\tau=\const\}
\end{equation}
are $C^{3,\al}$ spacelike hypersurfaces, or even smooth hypersurfaces, if $N$ is smooth, and $\tau$ is also continuous, since
\begin{equation}
\tau(x^0,x)=\fv H{M(t)}\qq\A\, (x^0,x)\in N,
\end{equation}
where
\begin{equation}
(x^0,x)\in M_t
\end{equation}
or, equivalently,
\begin{equation}\lae{8.5}
x^0=u(t,x),
\end{equation}
and we know that, if
\begin{equation}
(x^0_k,x_k)\ra (x^0_0,x_0)\qq\tup{in}\; N
\end{equation}
such that
\begin{equation}
(x^0_k,x_k)\in M_{t_k},
\end{equation}
then
\begin{equation}
t_k\ra t_0
\end{equation}
and
\begin{equation}
(x^0_0,x_0)\in M_{t_0},
\end{equation}
or equivalently,
\begin{equation}
x^0_0=u(t_0,x_0),
\end{equation}
in view of \frc{7.3}.

In order to prove that $\tau$ is of class $C^1$, it suffices to show 
\begin{equation}\lae{8.11}
u(\cdot,x)\in C^1(I)\qq\A\, x\in \so
\end{equation}
and
\begin{equation}\lae{8.12}
\dot u=\pde ut>0.
\end{equation}
Here, we already used the fact that the function in \re{8.11} is strictly monotone growing in $t$.

Indeed, if  \re{8.1} and \re{8.12} are valid, then we deduce from \re{8.5} 
\begin{equation}
\pde{x^0}t=\dot u>0.
\end{equation}
On the other hand, the function $\hat H$ in \re{8.1} is a diffeomorphism, i.e.,
\begin{equation}
\pde \tau t=\hat H '>0
\end{equation}
and hence we conclude
\begin{equation}
\pde{x^0}\tau=\pde{x^0}t \pde t\tau=\dot u (\hat H')^{-1}>0
\end{equation}
and
\begin{equation}
\pde\tau{x^0}=\frac{\hat H'}{\dot u}.
\end{equation}
Once we have proved \re{8.11} and \re{8.12} the conclusion $\tau$ is smooth if $N$ is smooth  follows immediately as we shall show.

Thus, let us prove \re{8.11} and \re{8.12} by showing that these relations are valid for the solutions $u=u_R$ of the Dirichlet problems \fre{3.1}. Since we shall  also derive 
\begin{equation}
0<\dot u_R(t,x)\le c\qq\A\, x\in B_R,
\end{equation}
where the constant $c$ is dependent of $R$ and also independent of $t$, als long as $t$ stays in a compact subset of $I$, we then conclude
\begin{equation}
0\le \dot u\le c\qq\A\, x\in \so,
\end{equation}
and the weak Harnack inequality will then allow us to actually deduce
\begin{equation}
0<\dot u.
\end{equation}
Thus, let $u=u_R(t,x)$ be a solution of \re{3.1} with boundary values $t$ and mean curvature $\hat H$. First, we need to prove the uniqueness of $u$.
\bl
The solution $u\in C^{3,\al}(\bar B_R)$ of the Dirichlet problem \re{3.1} is unique.
\el
\bp
Let $M_i=\graph u_i$, $i=1,2$, be two solutions of \re{3.1} and assume that there exists $x\in B_R$ such that
\begin{equation}
u_2(x)<u_1(x).
\end{equation}
Then, there exists a maximizing future directed timelike geodesic $\ga=(\ga^\al)$ from $\bar M_2$ to $\bar M_1$ and $p_i\in \bar M_i$ such that
\begin{equation}\lae{8.21}
0<d_0=d(M_2,M_1)=d(p_2,p_1),
\end{equation}
where
\begin{equation}
p_i=u_i(x_i),\qq\, x_i\in \bar B_R.
\end{equation}
Note that
\begin{equation}
d_0=\sup{\set{d(x_2,x_1)}{x_2\in M_2, x_1\in M_1}}
\end{equation}
and
\begin{equation}
0<d(p_2,p_1)\le d_0,
\end{equation}
where
\begin{equation}
p_i=(u_i(x),x),\qq i=1,2. 
\end{equation}
Since the $\bar M_i$ are compact and $d$ continuous, the supremum is achieved, hence \re{8.21} is valid and, furthermore, there exists a future directed timelike geodesic from $p_2$ to $p_1$, \cf \cite[Proposition 6.7.1, p. 212]{he:book}.

Let $Q$ be the open cylinder
\begin{equation}
Q=I\times B_R
\end{equation}
and let $\ga$ be parametrized in the interval $[0,1]$. Suppose
\begin{equation}
p_1\in \pa Q,
\end{equation}
then
\begin{equation}
p_2\not\in \pa Q,
\end{equation}
since $\ga$ is future directed and hence 
\begin{equation}
\ga^0(\tau)<\ga^0(1)=t\qq\A\, 0\le \tau<1.
\end{equation}
Let $\tau_0$ be the first
\begin{equation}
\tau\in (0,1]
\end{equation}
such that
\begin{equation}
\ga(\tau)\in \pa Q.
\end{equation}
If $\tau_0=1$, then
\begin{equation}
\ga(\tau)\in \bar Q\qq\A\,\tau\in [0,1].
\end{equation}
Since $M_2$ can be extended as a spacelike graph $\tilde M_2$ over a slightly larger ball $\bar B_{R+\e}$ $\ga$ would then be a future directed curve completely contained in the open cylinder
\begin{equation}
\tilde Q=I\times B_{R+\e}
\end{equation}
with endpoints $p_i\in \tilde M_2$, which leads to a contradiction, in view of \frl{3.2}.

Thus, let us assume
\begin{equation}
\tau_0<1,
\end{equation}
then
\begin{equation}
\ga^0(\tau_0)<u_2(\ga^i(\tau_0))=t
\end{equation}
and we can define a new future directed broken curve $\tilde \ga$ completely contained in $\tilde Q$ with endpoints in $\tilde M_2$, \cf part (ii) of the proof of \frl{3.1}, which again leads to a contradiction as before.

By switching the time orientation the previous arguments also exclude the case
\begin{equation}
p_2\in \pa Q.
\end{equation}
Thus, we may assume
\begin{equation}
p_i\in M_i,\qq i=1,2.
\end{equation}
But then we are in the same situation as in the proof of \rl{7.2} after equation \fre{7.39} and the arguments there  lead to a contradiction.
\ep

\bc\lac{8.2}
Let us write solutions $u_R$ of \re{3.1} in the form
\begin{equation}
u_R=u_R(t,x),\qq\, (t,x)\in I\times \bar B_R,
\end{equation}
then $u_R$ is continuous in $t$.
\ec
\bp
Obvious.
\ep

Next, let us prove:
\bl\lal{8.3}
The function $u_R(\cdot,x)$ belongs to $C^1(I)$ for all $x\in \bar B_R$.
\el
\bp
Let $t_0\in I$ be arbitrary, define
\begin{equation}
M_t=\graph u_R(t,\cdot)
\end{equation}
and set
\begin{equation}
M_0=M_{t_0}.
\end{equation}
The $\bar M_t$ are of class $C^{3,\al}$. Let $\mc U$ be a tubular neighbourhood of $M_0$ and $\vt$ the signed distance function to $M_0$ such that $\vt$ is a time function in $\mc U$; $\vt$ is also the time coordinate of a normal Gaussian coordinate system $(\vt,x^i)$ in $\mc U$, \cf \cite[Theorem 12.5.13]{cg:aII}. Since $M_0$ is a graph over $B_R$ we may consider the $(x^i)$ to be local coordinates in $B_R$. Note that
\begin{equation}
\vt\in C^{3,\al}(\mc U).
\end{equation}
Since we can extend $M_0$ to be a graph over a slightly larger ball we shall assume that $\mc U$ covers $\bar M_0$. Moreover, we know that there is $\de>0$ such that
\begin{equation}\lae{8.39} 
\bar M_t\su \mc U\qq\A\; t\in (t_0-\de,t_0+\de).
\end{equation}
In terms of $(x^0,x)$ $\vt$ can be expressed as
\begin{equation}\lae{8.40} 
\vt=\tilde\f(x^0,x),\qq \tilde\f\in C^{3,\al},
\end{equation}
and similarly
\begin{equation}\lae{8.41}
x^0=\f(\vt,x),\qq \f\in C^{3,\al}.
\end{equation}
Note that the spatial coordinates are the same in both coordinate systems. The slice
\begin{equation}\lae{8.42}
\{x^0=t\}\ii\mc U
\end{equation}
can then be expressed as the graph
\begin{equation}\lae{8.43} 
\{(\vt,x):\vt=\tilde \f(t,x)\}
\end{equation}
and the hypersurfaces $M_t$ in \re{8.39} as
\begin{equation}
\{(\vt,x):\vt=\tilde\f(u(t,x),x)\equiv \tilde u(t,x)\},
\end{equation}
and similarly
\begin{equation}\lae{8.45}
u(t,x)=\f(\tilde u(t,x),x).
\end{equation}
We shall now prove that
\begin{equation}
\tilde u(\cdot,x)\in C^1(I_{\de_0}),
\end{equation}
where
\begin{equation}
I_{\de_0}=(t_0-\de_0,t_0+\de_0)
\end{equation}
for some small
\begin{equation}
0<\de_0<\de.
\end{equation}
Let us define the function space
\begin{equation}
C^{3,\al}_0(\bar B_R)=\{u\in C^{3,\al}(\bar B_R):\fv u{\pa B_R}=0\}.
\end{equation}
In view of \re{8.42} and \re{8.43} the graph
\begin{equation}
\{\vt=\tilde\f(t,x)\}
\end{equation}
represents the graph
\begin{equation}
\{x^0=t\}
\end{equation}
in the new coordinate, hence
\begin{equation}
0=\tilde \f(t_0,x)\qq\A\, x\in \pa B_R
\end{equation}
and by continuity
\begin{equation}
-\e_1<\tilde\f(t,x)<\e_1\qq\A\, (t,x)\in \msc U_{\de_1},
\end{equation}
where
\begin{equation}
\msc U_{\de_1}=(t_0-\de_1,t_0+\de_1)\times \{x\in \bar B_R:R-\rho(x)<\de_1\}.
\end{equation}
Let
\begin{equation}
\chi\in C^\un_c(-\de_1,\de_1)
\end{equation}
be a real cut-off function such that
\begin{equation}
\chi(0)=1\q\wed\q 0\le \chi\le 1
\end{equation}
and set
\begin{equation}
w=\chi(R-\rho(x)),
\end{equation}
then
\begin{equation}
w(x) \tilde\f(\cdot, x)\in C^1(I_{\de_1}) \qq\A\, x\in \bar B_,
\end{equation}
where
\begin{equation}
I_{\de_1}=(t_0-\de_1,t_0+\de_1),
\end{equation}
and
\begin{equation}
R(w\tilde\f)\su J_\de=(-\de,\de)\qq\A\, (t,x)\in I_{\de_1}\times \bar B_R.
\end{equation}
Hence, the last relation is also valid for all
\begin{equation}
\tilde u+w\tilde\f(t,\cdot)\qq\A\, t\in I_{\de_1},
\end{equation}
where
\begin{equation}\lae{8.62} 
\tilde u\in \msc B_{\e_0}(-w\tilde\f)\su C^{3,\al}_0(\bar B_R),
\end{equation}
if $\e_0$ and $\de_1$ are small enough, and
\begin{equation}
M=\graph (\tilde u+w\tilde\f(t,\cdot))
\end{equation}
is a spacelike graph in $(\vt,x)$ which is contained in the tubular neighbourhood $\mc U$ which can be expressed in the new coordinate system as a cylinder
\begin{equation}
\mc U=(-\de,\de)\times B_{R+\e}.
\end{equation}
If the hypersurface $M$ is viewed as a graph in $(x^0,x)$,
\begin{equation}
M=\graph u,
\end{equation}
then
\begin{equation}
\fv u{\pa B_R}=t.
\end{equation}
Indeed, according to \re{8.45}
\begin{equation}
u(x)=\f(\tilde u(x)+w\tilde\f(t,x),x)
\end{equation}
and
\begin{equation}
u(x)=\f(\tilde\f(t,x),x)=t\qq\A\, x\in \pa B_R.
\end{equation}
We have
\begin{equation}
M_0=\{(\vt,x):\vt=0,\;x\in B_R\},
\end{equation}
i.e.,
\begin{equation}
M_0=\graph (-w \tilde\f(t_0,\cdot)+w \tilde\f(t_0,\cdot)).
\end{equation}
Hence, for all $\tilde u$ in \re{8.62} the hypersurfaces
\begin{equation}
M(t,\tilde u)=\graph(\tilde u+w\tilde\f(t,\cdot)),\qq t\in I_{\de_1},
\end{equation}
are spacelike and contained in $\mc U$. Thus, we can define the operator
\begin{equation}
\begin{aligned}
G:&I_{\de_1}\times \msc B_{\e_0}(-w\f(t_0,\cdot))\ra C^{1,\al}(\bar B_r),\\
&G(t,\tilde u)=\fv H{M(t,\tilde u)}-\hat H(t)
\end{aligned}
\end{equation}
and consider the equation
\begin{equation}
G(t,\tilde u)=0
\end{equation}
which describes the Dirichlet problem \fre{3.1} in the coordinate system $(\vt,x)$ as an implicit function equation. We know that $G$ is at least of class $C^1$ in $(t,\tilde u)$ and that
\begin{equation}
G(t_0,-w\tilde\f(t_0,\cdot))=0
\end{equation}
and
\begin{equation}
D_2G(t_0,-w\tilde\f(t_0,\cdot))=L,
\end{equation}
where $L$,
\begin{equation}
L:C^{3,\al}_0(\bar B_R)\ra C^{1,\al}(\bar B_R),
\end{equation}
is the elliptic differential operator
\begin{equation}\lae{8.77}
L\tilde u=-\D \tilde u+(\abs A^2+\bar R_{\al\bet}\nu^\al\nu^\bet)\tilde u.
\end{equation}
Here, the Laplacian is defined by the metric $g_{ij}$ in $M_0$,
\begin{equation}
\abs A^2=h^{ij}h_{ij},
\end{equation}
and $(\nu^\al)$ is the normal of $M_0$.

Since 
\begin{equation}
N(L)=\{0\}
\end{equation}
$L$ is a topological homeomorphism, and hence, we conclude from the implicit function theorem that there exists a small neighbourhood
\begin{equation}
I_\e(t_0)\su I_{\de_1}
\end{equation}
and a uniquely determined function
\begin{equation}
\theta\in C^1(I_\e(t_0),C^{3,\al}_0(\bar B_R))
\end{equation}
such that
\begin{equation}\lae{8.82}
G(t,\theta (t))=0\qq\A\, \in I_\e(t_0).
\end{equation}
Writing
\begin{equation}\lae{8.83} 
\tilde u(t,x)=\theta(t) +w\tilde\f(t,\cdot)
\end{equation}
we conclude that $\tilde u(\cdot,x)$ is of class $C^1$ in the interval $I_\e(t_0)$. Expressing the graphs in the original coordinates we obtain, in view of \re{8.45},
\begin{equation}\lae{8.84}
u(\cdot,x)=\f(\tilde u(\cdot,x),x)
\end{equation}
is of class $C^1$ in the intervall $I_\e(t_0)$.
\ep

\bc\lac{8.4}
From \re{8.84} we conclude
\begin{equation}
\dot u(t_0,x)=\dot\f(0,x)\dot{\tilde u}(t_0,x)=e^{-\psi}\tilde v \dot{\tilde u}(t_0,x)\qq\A
\, x\in \bar B_R.
\end{equation}
\ec
\bp
We only have to prove that
\begin{equation}
\dot\f(0,x)=e^{-\psi}\tilde v.
\end{equation}
Let us recall that the tubular neighbourhood is defined by the geodesic flow $\ga(\vt,x)=(\ga^\al)$ defined in the coordinates $(x^0,x^i)$. Here $\vt$ is the signed arc length which is identical to the signed distance function of $M_0$, the initial values of the flow at $\vt=0$ are
\begin{equation}
\ga(0,x)=(u(t_0,x),x)\q\wed\q \dot\ga (0,x)=-\nu(x),
\end{equation}
where $\nu$ is the past directed normal of $M_0$. Thus, the flow is future directed. The function $\f(\vt,x)$ in the relation \re{8.41} is then identical to
\begin{equation}
\f(\vt,x)=\ga^0(\vt,x),
\end{equation}
hence
\begin{equation}\lae{8.89}
\dot\f(0,x)=\dot\ga^0(0,x)=-\nu^0=e^{-\psi}\tilde v.
\end{equation}
\ep

\bpp\lap{8.5}
Let $\tilde u(t,x)=\tilde u_R(t,x)$ be a solution of the Dirichlet problem \re{3.1} in the tubular neighbourhood $\mc U$ of
\begin{equation}
M_0=\graph u_R(t_0,\cdot),
\end{equation}
then
\begin{equation}
0\le \dot{\tilde u}(t_0,x)\le c\qq\A\, x\in B_R,
\end{equation}
where $c$ depends on known estimates but not on $R$.
\epp
\bp
Let us differentiate the equations \re{8.82} and \re{8.83} with respect to $t$ and evaluate at $t=t_0$. Then we obtain
\begin{equation}\lae{8.92}
-\D\dot {\tilde u}(t_0) +(\abs A^2+\bar R_{\al\bet}\nu^\al\nu^\bet)\dot{\tilde u}(t_0)=\hat H'(t_0)
\end{equation}
as well as the boundary condition
\begin{equation}\lae{8.93}
\fv{\dot{\tilde u}(t_0,x)}{\pa B_R}=\fv{\dot{\tilde\f}(t_0,x)}{\pa B_R},
\end{equation}
where
\begin{equation}
\hat H'=\pde{\hat H}t>0.
\end{equation}
Let us first note that
\begin{equation}
t=\f(\tilde\f(t,x),x)
\end{equation}
and therefore
\begin{equation}
1=\dot\f(0,x)\dot{\tilde\f}(t_0,x),
\end{equation}
hence,
\begin{equation}
\dt{\f}(0,x)=e^\psi v,
\end{equation}
in view of \re{8.89}, i.e., the boundary values of $\dt{u}(t_0,x)$ are uniformly bounded and tend to $1$ if $R$ tends to infinity.

Secondly, the coefficient in \re{8.92} is non-negative
\begin{equation}
\abs A^2+\bar R_{\al\bet}\nu^\al\nu^\bet\ge 0
\end{equation}
and strictly positive for large $r(x)$, or equivalently, $\rho(x)$,
\begin{equation}
\abs A^2+\bar R_{\al\bet}\nu^\al\nu^\bet\ge c_0>0\qq\A\, x\in\{\rho(x)\ge R_1\},
\end{equation}
where $R_1$ is sufficiently large , \cf \frc{6.4}.

Let us now estimate $\dt u=\dt u(t_0,\cdot)$ from above. Assuming $R>R_1$ we choose $k_0\in \R[]$ such that
\begin{equation}\lae{8.100}
k_0\ge \sup_{\pa B_R}e^\psi v
\end{equation}
and define for $k\ge k_0$
\begin{equation}\lae{8.101}
\h=\max(\dt u-k,0)\in H^{1,2}_0(B_R).
\end{equation}
Multiplying \re{8.92} by $\h$ and integrating by parts we obtain
\begin{equation}
\int_{B_R}\abs {D\h}^2+\int_{B_R}(\abs A^2+\bar R_{\al\bet}\nu^\al\nu^\bet)\dt u\h=\int_{B_R}\hat H'\h.
\end{equation}
Choosing at the moment $k=k_0$ and assuming furthermore that $k_0$ satisfies besides \re{8.100} also
\begin{equation}\lae{8.103} 
c_0k_0\ge \hat H'(t_0)+1
\end{equation}
we infer
\begin{equation}\lae{8.104}
\int_{B_R}\abs {D\h}^2+\int_{B_R\sminus B_{R_1}}\h\le \int_{B_{R_1}}\hat H'\h.
\end{equation}
The right-hand side can be estimated by 
\begin{equation}
\int_{B_{R_1}}\hat H'\h\le \hat H'\abs{B_{R_1}}^\frac {2n}{n+2}\Big(\int_{B_R}\abs \h^\frac{2n}{n-2}\Big)^\frac{n-2}{2n}.
\end{equation}
On the other hand, in view of the Sobolev embedding theorem, we have
\begin{equation}
\Big(\int_{B_R}\abs \h^\frac{2n}{n-2}\Big)^\frac{n-2}{n}\le c \int_{B_R}\abs{D\h}^2
\end{equation}
and thus we conclude
\begin{equation}
\int_{B_R}\abs{D\h}^2\le c=\const
\end{equation}
independent of $R$. This inequality is also valid if $\h$ is defined as in \re{8.101} for $k\ge k_0$.

Consider now $k\ge k_0$ to be arbitrary and let us look again at the inequality \re{8.104}. Define
\begin{equation}
A(k)=\{x\in B_R:\dt u>k\},
\end{equation}
then we deduce
\begin{equation}
\Big(\int_{B_R}\abs \h^\frac{2n}{n-2}\Big)^\frac{n-2}{n}\le c \abs{A(k)}^\frac{n+2}{2n}\Big(\int_{B_R}\abs \h^\frac{2n}{n-2}\Big)^\frac{n-2}{2n}
\end{equation}
implying
\begin{equation}\lae{8.110}
\Big(\int_{B_R}\abs \h^\frac{2n}{n-2}\Big)^\frac{n-2}{2n}\le c \abs{A(k)}^\frac{n+2}{2n}.
\end{equation}
Next, let $h>k$ then
\begin{equation}
(h-k)\abs{A(h)}\le \int_{B_R}\h\le c \abs{A(k)}^{1+\frac2n}\qq\A\, h>k\ge k_0
\end{equation}
from which we conclude, in view of a lemma due to Stampacchia, \cf \cite[Lemma 4.1, p. 93]{stampacchia:montreal},
\begin{equation}
\dt u\le k_0+d,
\end{equation}
where
\begin{equation}
d=c\abs{A(k_0)}^\frac2n 2^{\frac n2+1}.
\end{equation}
To prove the uniform boundedness of $\abs{A(k_0)}$ we assume that $k_0$ is so large that
\begin{equation}
k_1=\frac{k_0}2-1
\end{equation}
also satisfies the inequalities in \re{8.100} and \re{8.103} if we replace $k_0$ by $k_1$. Then
\begin{equation}
\int_{B_R}\h_0^\frac{2n}{n-2}\le c,
\end{equation}
where
\begin{equation}
\h_0=\max(\dt u-k_1,0),
\end{equation}
and we deduce
\begin{equation}
\begin{aligned}
\abs{A(k_0)}&\le \frac1{(\frac{k_0}2+1)^\frac{2n}{n-2}}\int_{A(k_0)}(u-k_1)^\frac{2n}{n-2}\\
&\le \frac1{(\frac{k_0}2+1)^\frac{2n}{n-2}}\int_{B_R}\abs{\h_0}^\frac{2n}{n-2}.
\end{aligned}
\end{equation}
\ep

\bpp\lap{8.6}
Let $\dt u(t_0,x)$ be the solution of \re{8.92} and the corresponding boundary condition, then
\begin{equation}
\abs{\dt u(t_0,\cdot)}_{2,\al,\bar B_R}\le c,
\end{equation}
 where $c$ only depends on already proven estimates but not on $R$.
\epp
\bp
The proof is identical to the proof of \frl{6.2}.
\ep

If we want to improve the regularity of $\tilde u(t,x)$ we have to assume a higher regularity of $N$. At the moment we only assume $N$ to be of class $C^{3,\al}$. Thus, let us suppose $N$ to be of class $C^{m,\al}$, $m\ge 3$, $0<\al<1$. Then, the original coordinates $(x^\al)$ with the time function $x^0$ are of class $C^{m,\al}$, the metric $\bar g_{\al\bet}$ of class $C^{m-1,\al}$ and the second fundamental form $\bar h_{ij}$ of the coordinate slices and the Riemann curvature tensor are of class $C^{m-2,\al}$ with uniform bounds. 

The solutions $u(t,x)=u_R(t,x)$ of the Dirichlet problems \fre{3.1} are then of class $C^{m,\al}(\bar B_R)$ with uniform bounds independent of $R$. The estimates, of course, depend on $m$.

First, let us prove:
\bl\lal{8.7}
Let $N$ be of class $C^{m,\al}$, $m\ge 3$, $0<\al<1$, with uniform bounds, $u(t,x)=u_R(t,x)$ a solution of the Dirichlet problem \re{3.1} with boundary value $t_0$ and let 
\begin{equation}
M_0=\graph u(t_0,\cdot).
\end{equation}
Let $(\vt,x)$ be the normal Gaussian coordinate system corresponding to a tubular neighbourhood $\mc U$ of $\bar M_0$, then $\vt$ and the transformation maps
\begin{equation}
x^0=\f(\vt,x)
\end{equation}
and its inverse
\begin{equation}
\vt=\tilde\f(x^0,x)
\end{equation}
are of class $C^{m,\al}$ with respect to the indicated variables such that the corresponding $C^{m,\al}$-norms, evaluated at $\bar M_0$, are uniformly bounded independent of $R$.
\el

\bp
It suffices to prove the claim for $\f$.  As we mention before, the new coordinate system is created with the help of the geodesic flow $\ga(\vt,x)$ with initial hypersurface $\bar M_0$ and initial values
\begin{equation}\lae{8.122}
\ga(0,x)=(u(x),x)\q\wed\q \dot\ga(0,x)=-\nu(x).
\end{equation}
$\vt$ is the signed arc length, which will be negative in the past of $M_0$, $\nu$ is supposed to be  the past directed normal of $M_0$. The map $\f$ is defined by
\begin{equation}
x^0=\ga^0(\vt,x)\equiv \f(\vt,x).
\end{equation}
In view of \re{8.122} we immediately obtain 
\begin{equation}
\ga(0,\cdot)\in C^{m,\al}(\bar B_R)\q\wed\q \dot \ga (0,\cdot)\in C^{m-1,\al}(\bar B_R)
\end{equation}
with uniform bounds independent of $R$ and from the geodesic equation and  the assumption on $N$ we recursively deduce
\begin{equation}
D^k_\vt\ga (0,\cdot)\in C^{m-k,\al}(\bar B_R)\qq\A\, 0\le k\le m,
\end{equation}
again with uniform bounds independent of $R$.
\ep

We can now prove:
\bt\lat{8.8}
Let $N$ be of class $C^{m,\al}$, $m\ge 3$, $0<\al<1$, and let $\tilde u(t,x)=\tilde u_R(t,x)$ be the solutions of the Dirichlet problems \re{3.1} expressed in the new coordinates $(\vt,x)$. Then $\tilde u$ is of class $C^{m-2}$ with respect to $t$ and
\begin{equation}
\frac{\pa^k\tilde u}{(\pa t)^k}(t_0,\cdot)\in C^{m-k,\al}(\bar B_R)\qq\A\, 0\le k\le m-2
\end{equation}
with uniform bounds independent of $R$.
\et
\bp
We shall prove the theorem inductively; for $k=0,1$ this has already been proved before, \cf \rp{8.6} and apply the Schauder estimates to the solution of \re{8.92} with boundary values given in \re{8.93}. 

To prove the claim for $k=2$, let  us consider the elliptic differential equation satisfied by $\tilde u(t,\cdot)$ for $t$ near $t_0$
\begin{equation}\lae{8.127}
\begin{aligned}
\hat H(t) \tilde v&=-g^{ij}\tilde u_{ij}+g^{ij}\bar h_{ij}\\
&=-g^{ij}\tilde u_{ij}+\bar H+\tilde v^2\bar h_{ij}\check{\tilde u}^i\check{\tilde u}^j.
\end{aligned}
\end{equation}
Here, $g_{ij}$ is the induced metric of $M_t=\graph \tilde u(t,\cdot)$ and the ambient metric is expressed in normal Gaussian coordinates
\begin{equation}
d\bar s^2=-d\vt^2+\s_{ij}(\vt,x)dx^idx^j
\end{equation}
in the tubular neighbourhood $\mc U$ of $M_0$. At $t=t_0$ we know
\begin{equation}
\tilde u(t_0,\cdot)\equiv 0.
\end{equation}
Differentiating \re{8.127} twice with respect to $t$, the differentiability is due to the implicit function theorem, and evaluating the result at $t=t_0$ we obtain
\begin{equation}
\hat H''+F=-g^{ij}\Ddot {\tilde u}_{ij}+\dot{\bar H}\Ddot{\tilde u},
\end{equation}
where
\begin{equation}
F=F(x,\tilde u, \dt u, D\tilde u,  D\dt u, D^2\tilde u, D^2\dt u)\in C^{m-3,\al}(\bar B_R)
\end{equation}
with uniform bounds independent of $R$. The proof of \rp{8.5} then yields
\begin{equation}
\abs {\Ddot{\tilde u}(t_0,\cdot)}\le c \qq\tup{in}\; \bar B_R,
\end{equation}
where, now, we also have to estimate $\Ddot{\tilde u}$ from below. Combining this estimate with the Schauder estimates we conclude
\begin{equation}
\abs{\Ddot{\tilde u}(t_0,\cdot)}_{m-2,\al,\bar B_R}\le c,
\end{equation}
since
\begin{equation}
\fv{\Ddot{\tilde u}}{\pa B_R}\in C^{m-2,\al}.
\end{equation}

The same arguments also apply when higher derivatives are considered. Set
\begin{equation}
\tilde u^{(k)}=\frac{\pa^k\tilde u}{(\pa t)^k}(t_0,\cdot)
\end{equation}
then
\begin{equation}
D^k_t\hat H+F=-\D\tilde u^{(k)}+\dot{\bar H}\tilde u^{(k)},
\end{equation}
where, now, $F$ depends on
\begin{equation}
F=F(x,\tilde u,\dot{\tilde u}, \ldots, D^2\tilde u^{(k-1)})\in C^{m-k,\al}(\bar B_R),
\end{equation}
and the boundary values are of class $C^{m-k,\al}$.
\ep

\bc \lac{8.9}
The results of \rt{8.8} are also valid for $u_R(t_0,x)$, in view of \rl{8.7}, see also \rc{8.4}. Moreover, since $t_0$ is arbitrary, these estimates are valid for any $t$ and they are uniform provided $t$ ranges in a compact subset of $I$.
\ec

Letting $R$ tend to infinity we then deduce:
\bt\lat{8.10}
The functions
\begin{equation}
u(t,x),\qq (t,x)\in I\times \so,
\end{equation}
describing the foliation hypersurfaces $M_t$, $t\in I$, are of class $C^{m-3,1}$ in $t$ such that
\begin{equation}
D^k_t u\in C^{m-k,\al}(\so)\qq\A\, 0\le k\le m-3,
\end{equation}
if $N$ is of class $C^{m,\al}$, $m\ge 3$, $0<\al<1$; if $N$ is smooth, i.e., $m=\un$, then $u$ is also smooth in the variables $(t,x)$ and the mean curvature function $\tau=\tau(x^0,x)$ is a smooth time function.
\et
\bp
This immediately follows from \frl{6.2} and the considerations at the beginning of this section.
\ep

\section{Appendix: Lipschitz continuous solutions are regular}\las{9}

In this appendix we want to prove that Lipschitz continuous solutions of the Dirichlet problems \fre{3.1} or of the related equation \fre{5.96} are of class $H^{2,p}$ or even more regular depending on the data. 

We shall use the assumptions stated at the beginning of the proof of \rl{5.4} and consider an allowed Lipschitz continuous solution $u$ of \re{5.96} in $B_R$ with vanishing boundary values which satisfies
\begin{equation}
\abs{Du}^2=\s^{ij}(u,x)D_iuD_ju\le 1-\de^2,
\end{equation}
and
\begin{equation}
m_1\le u\le m_2,\qq m_i\in I.
\end{equation}
The equation \re{5.96} has the form
\begin{equation}\lae{9.3}
-D_i(a^i(x,u,Du))=f^0,
\end{equation}
where the divergence is with respect to the metric
\begin{equation}
\s_{ij}(u,x).
\end{equation}
Let $x^0$ be a fixed constant, then we can express the divergence in \re{9.3} with respect to metric
\begin{equation}
\s_{ij}(x^0,x)
\end{equation}
without changing the structure or the properties of the coefficients and $f^0$. The volume element in the integrations below will also be defined by this metric.

We now argue similarly as in \cite[Section 1]{cg85} and modify the coefficients of the operator---only a slight adaptation to the present situation is necessary. First, let $\vt$ be a smooth real function such that
\begin{equation}
\vt(t)=
\begin{cases}
\msp[83]t, &m_1\le t\le m_2,\\
-(m_1+\e_1),&t\le -(m+\e_1),\\
\hp{-}(m_2+\e_1),&t\ge (m_2+\e_1),
\end{cases}
\end{equation}
where $\e_1>0$ is chosen small enough  to guarantee that
\begin{equation}
\vt(\R[])\su I.
\end{equation}
Secondly, define the metric $\tilde \s_{ij}(t,x)$ by
\begin{equation}
\tilde\s_{ij}(t,x)=\s_{ij}(\vt(t),x).
\end{equation}
Thirdly, if a vector field $p=(p_i)$ is the gradient of a function $w$,
\begin{equation}
p=Dw,
\end{equation}
then set
\begin{equation}
\abs p^2=\tilde\s^{ij}p_ip_j=\tilde\s^{ij}(w,x)D_iwD_jw.
\end{equation}

Let $\om$, $g$ be smooth real functions such that
\begin{equation}
\om(t)=
\begin{cases}
1,&0\le t\le 1-\frac{\de^2}2,\\
0,&t\ge 1-\frac{\de^2}3,
\end{cases}
\end{equation}
and assume $g$ to be convex satisfying
\begin{equation}
g(t)=
\begin{cases}
\hp{c}0,&0\le t\le 1-\de^2,\\
ct ,&t\ge 1-\frac{\de^2}2,
\end{cases}
\end{equation}
where $c$ is some positive constant. Then, we define 
\begin{equation}
\tilde a^i(x,t,p)=a^i(x,\vt(t),p)\om(\abs p^2)+kg'(\abs p^2)p^i.
\end{equation}
Here, $k$ is a positive constant. It can  be easily verified that
\begin{equation}
\tilde a^{ij}=\frac{\pa\tilde a^i}{\pa p_j}
\end{equation}
is uniformly elliptic if $k$ is large enough.

Then, we look at the Dirichlet problem
\begin{equation}\lae{9.12}
\begin{aligned}
\tilde Aw=-D_i(\tilde a^i(Dw))+\ga (w-u)&=f^0,\\
\fv w{\pa B_R}&=0,
\end{aligned}
\end{equation}
where $\ga>0$ is so large that the operator on the left-hand side of \re{9.12} is uniformly monotone, i.e., 
\begin{equation}
\spd{\tilde Aw_2-\tilde Aw_1}{w_2-w_1}\ge c_0\norm{w_1-w_2}^2\qq\A\, w_1,w_2\in H^{1,2}_0(B_R),
\end{equation}
where $c_0$ is positive. The pairing is the bilinear form between $H^{1,2}_0(B_R)$ and its dual space $H^{-1,2}(B_R)$. 

Evidently, the solutions of the Dirichlet problem are uniquely determined and since
\begin{equation}
\tilde Au=-D_i(a^i(x,u,Du))=f^0,
\end{equation}
and
\begin{equation}\lae{9.19.1}
\tilde \s_{ij}(u,x)=\s_{ij}(u,x),
\end{equation}
we deduce that a solution of \re{9.12} has to coincide with $u$. It is well known, due to the Calderon-Zygmund inequalities and the De Giorgi-Nash theorem,  that \re{9.12} has a solution
\begin{equation}\lae{9.19}
w\in H^{2,p}(B_R)
\end{equation}
for any $1<p<\un$. Indeed,  since the operator $\tilde A$ is uniformly elliptic and monotone, the unique solution $w\in H^{1,2}_0(B_R)$ is then also of class $H^{2,2}(B_R)$, because of the $L^2$-estimates, and, since $w=u$, $w$ is Lipschitz. Furthermore, $Dw$ is H\"older continuous with exponent $\bet$ for some $0<\bet<1$ with uniform H\"older norm in
\begin{equation}
B_{\rho_0}(x_0)\su \bar B_{2\rho_0}(x_0)\su B_R
\end{equation}
and also in
\begin{equation}
\Om(x_0,\rho_0),\qq\, x_0\in\pa B_R,
\end{equation}
in view of the De Giorgi-Nash theorem. Hence, the coefficients are continuous and the Calderon-Zygmund inequalities finally yield \re{9.19}. If the data are better we can then  apply the Schauder estimates.

\br\lar{9.1}
The norm of $u$ in $H^{2,p}(B_R)$ depends on $R$, but local $H^{2,p}$ norms in
\begin{equation}
B_{\rho_0}(x_0)\su \bar B_{2\rho_0}(x_0)\su B_R
\end{equation}
and also in
\begin{equation}
\Om(x_0,\rho_0),\qq\, x_0\in\pa B_R,
\end{equation}
only depend on $\rho_0$ and $p$ but not on $R$, hence, for any $0<\al<1$
\begin{equation}
\abs u_{1,\al, B_R}\le c
\end{equation}
uniformly in $R$.
\er

%\backmatter
%\includepdf[pages=-]{/Users/claus/Documents/Scanned-Documents/}
%\bibliographystyle{hamsplain}
%\bibliography{mrabbrev,publications}
%\bibliography{publications}
\providecommand{\bysame}{\leavevmode\hbox to3em{\hrulefill}\thinspace}
\providecommand{\href}[2]{#2}

%\listoffigures

%\cleardoublepage

%\thispagestyle{empty}
%\closegraphsfile
\end{document}